\documentclass[preprint,12pt]{elsarticle}
\usepackage[latin1]{inputenc}
\usepackage{amsmath}
\usepackage{graphics}
\usepackage{xcolor}
\usepackage{amssymb}
\usepackage{amsmath, amsthm, amssymb}
\usepackage{url}
\graphicspath{{fig/}}     

\newtheorem{lem}{Lemma}[section]

\newtheorem{thm}{Theorem}
\newtheorem{lemma*}{Lemma}

\journal{Systems and Control Letters}
\makeatletter

\usepackage{pslatex,amssymb,amsthm}

\begin{document}

\begin{frontmatter}

\title{Output Feedback Control of Radially-Dependent Reaction-Diffusion PDEs on Balls of Arbitrary Dimensions\tnoteref{mytitlenote}}
\tnotetext[mytitlenote]{This work was supported in part by the National Natural Science Foundation of China (61773112), the Fundamental Research Funds for the Central Universities and Graduate Student Innovation Fund of Donghua University (CUSF-DH-D-2019089) and the scholarship from China Scholarship Council (CSC201806630010). Rafael Vazquez acknowledges financial support of the Spanish Ministerio de Ciencia, Innovaci\'on y Universidades under grant PGC2018-100680-B-C21}

\author{Rafael Vazquez\corref{corr}}\cortext[corr]{Corresponding author. Phone number +34954488148.}

\address{Dept. of Aerospace Engineering, Universidad de Sevilla, Camino de los Descubrimiento s.n., 41092 Sevilla, Spain.}

\ead{rvazquez1@us.es}

\author{Jing Zhang, Jie Qi\corref{jzjq}}

\address{College of Information Science and Technology, Donghua University, Shanghai 201620 China}

\ead{zhangj0811@163.com,jieqi@dhu.edu.cn}

\author{Miroslav Krstic\corref{mk}}

\address{Dept. of Mechanical and Aerospace Engineering, University of California San Diego, La Jolla, CA 92093-0411, USA.}

\ead{krstic@ucsd.edu}
%

\begin{abstract}
Recently, the problem of boundary stabilization and estimation for unstable linear constant-coefficient reaction-diffusion equation on $n$-balls (in particular, disks and spheres) has been solved by means of the backstepping method. However, the extension of this result to spatially-varying coefficients is far from trivial. Some early success has been achieved under simplifying conditions, such as radially-varying reaction coefficients under revolution symmetry, on a disk or a sphere. These particular cases notwithstanding, the problem remains open. The main issue is that the equations become singular in the radius; when applying the backstepping method, the same type of singularity appears in the kernel equations. Traditionally, well-posedness of these equations has been proved by transforming them into integral equations and then applying the method of successive approximations. In this case, with the resulting integral equation becoming singular, successive approximations do not easily apply. This paper takes a different route and directly addresses the kernel equations via a power series approach, finding in the process the required conditions for the radially-varying reaction (namely, analyticity and evenness) and showing the existence and convergence of the series solution. This approach provides a direct numerical method that can be readily applied, despite singularities, to both  control and observer boundary design problems.
\end{abstract}


\end{frontmatter}

 \interdisplaylinepenalty=2500

\section{Introduction}
In this paper we introduce an explicit  boundary output-feedback control law to stabilize an unstable linear \emph{radially-dependent} reaction-diffusion equation on an $n$-ball (which in 2-D is a disk and in 3-D a sphere).

This paper extends the spherical harmonics~\cite{atkinson} approach of~\cite{nball}, which assumed constant coefficients, using some of the ideas of~\cite{vazquez-zhang}; for the sake of brevity we will mainly show the modifications required with respect to~\cite{nball}, skipping the details when they are identical. For a finite number of harmonics, we design boundary feedback laws and output injection gains using the backstepping method~\cite{krstic} (with kernels computed using a power series approach) which allows us to obtain exponential stability of the origin in the $L^2$ norm. Higher harmonics will be naturally open-loop stable. The required conditions for the radially-varying coefficients are found in the analysis of the numerical method and are non-obvious (evenness of the reaction coefficient). The idea of using a power series to compute backstepping kernels was first seen in~\cite{ascencio} (without much analysis of the method itself, but rather numerically optimizing the approximation) and later in~\cite{leo}, where piecewise-smooth kernels require the use of several series. Here, we prove that the method provides a unique converging solution.

Some partial results towards the solution of this problem were obtained in~\cite{catalan} and~\cite{sphere} for the disk and sphere, respectively; however they required simmetry conditions. Older results in this spirit were obtained in~\cite{Prieur} and~\cite{scott}. This paper extends and completes our conference contribution~\cite{jing-paper} where the ideas where initially presented (without proof).

Previous results and applications in multi-dimensional domains include  multi-agent deployment in 3-D space~\cite{jie} (by combining the ideas of~\cite{nball} and~\cite{meurer}), convection problem on annular domains~\cite{convloop}, PDEs with boundary conditions governed by lower-dimensional PDEs~\cite{jie2, vazquez-zhang}, multi-dimensional cuboid domains~\cite{meurer2}.

The backstepping method has proved itself to be an ubiquitous method for PDE control, with many other applications including, among others, flow control~\cite{vazquez,vazquez-coron}, nonlinear PDEs~\cite{vazquez2}, hyperbolic 1-D systems~\cite{long-nonlinear,auriol,rijke}, or delays~\cite{krstic5}. Nevertheless, other design methods are also applicable to the geometry considered in this paper (see for instance~\cite{triggiani} or~\cite{Barbu}).

The structure of the paper is as follows. In Section~\ref{sec-plant} we introduce the problem. In Section~\ref{sec-stab} we state our stability result. We study the well-posedness of the kernels in Section~\ref{sec-wp}, which is the main result of the paper, proving existence of the kernels and providing means for their computation; interestingly, odd and even dimensions require a slightly different approach. We briefly talk next about the observer in Section~\ref{sec-observer}, but skip most details based on its duality with respect to the controller. Then, we give some simulation results in Section~\ref{simulation}. We finally conclude the paper with some remarks in Section~\ref{sec-conclusions}.

\section{$n$-D Reaction-Diffusion System on an $n$-ball}\label{sec-plant}
Following~\cite{nball}, a varying coefficient reaction-diffusion system in an $n$-dimensional ball of radius $R$ can be written in $n$-dimensional spherical coordinates,
also known as ultraspherical coordinates (see~\cite{atkinson}, p. 93), which consist of one radial coordinate and $n-1$ angular coordinates. Then, using  a (complex-valued) Fourier-Laplace series of Spherical Harmonics\footnote{Spherical harmonics were introduced by Laplace to solve the homonymous equation and have been widely used since, particularly in geodesics, electromagnetism and computer graphics. A very complete treatment on the subject can be found in~\cite{atkinson}.} to handle the angular dependencies, one reaches the following independent complex-valued 1-D reaction-diffusion equation for each harmonic:
\begin{equation}\label{eqn-un}
\partial_t u^{m}_{l}=\frac{\epsilon}{r^{n-1}} \partial_r \left( r^{n-1} \partial_r u^{m}_{l} \right)-l(l+n-2) \frac{\epsilon}{r^2} u^{m}_{l}+\lambda(r) u^{m}_{l},
\end{equation}
evolving in $r \in[0,R],\,t>0$, with  boundary conditions 
\begin{eqnarray}\label{eqn-un-bc}
u^{m}_{l}(t,R)&=&U^{m}_{l}(t),
\end{eqnarray}
In these equations, we have considered Dirichlet boundary conditions. The measurement would be the flux at the boundary, namely  $\partial_r u_l^m (t,R)$.

In the above equations, the integers $m$ and $l$ stand for the order and degree of the harmonic, respectively. Note that the higher the degree  (corresponding to high
frequencies), the more ``naturally'' stable (\ref{eqn-un})--(\ref{eqn-un-bc}) is, as seen next. Define the $L^2$ norm  
\begin{equation}
	\rVert f(r) \rVert_{L^2}^2=\int\limits_{0}^{R}\rvert f(r)\rvert^2 r^{n-1}dr.
\end{equation}
and the associated $L^2$ space as usual, where $\vert f \vert^2 =ff^* $,
being $f^*$ the complex conjugate of $f$.

\begin{lem}\label{lem_ul_open}
	Given $\lambda(r)$ and $R$, there exists $L \in \mathbb{N}$ such that, for all $ l  >L$, the equilibrium $u_l^m\equiv 0$ of system \eqref{eqn-un}-\eqref{eqn-un-bc} is open loop exponentially stable, namely, for $U_l^m=0$ in \eqref{eqn-un} there exists a positive constant  $D_1$, such that for all $t$
	\begin{align}
		\label{equivalent2}\lVert u_l^m(t,\cdot) \lVert^2_{L^2} &\leq e^{-D_1t}  \lVert u_l^m(0,\cdot) \lVert^2_{L^2}.
	\end{align}
$D_1$ is independent of $l$, and only depends on $n$, $\lambda(r)$, $\epsilon$, and $R$, and can be chosen as large as desired just by increasing the values of $L$.
\end{lem}
The proof is skipped as it mimics~\cite{vazquez-zhang} just by using the $L^2$ norm as a Lyapunov function and Poincare's inequality.

Thus one only needs to stabilize the unstable mode with $l  < L$. Since the different modes are not coupled, it allows us to stabilize them separately and re-assembling them. Moreover since only a finite number of harmonics is stabilized, there is no need to worry about the convergence of the control law as in~\cite{nball}, with its Spherical Harmonics series being just a finite sum.

Our objective can now be stated as follows. Considering only the unstable modes, design an \emph{output-feedback} control law for $U_l^m$ using, for each mode, only the measurement of $\partial_r u_l^m (t,R)$. Our design procedure is established in the next section along with our main stability result.

\section{Stability of controlled harmonics}\label{sec-stab}

Next, for the unstable modes we design the output-feedback law. The observer and controller are designed separately using the backstepping method, by following~\cite{nball}; in this reference it is shown that both the feedback and the output injection gains can be found by solving a certain kernel PDE equation, which is essentially the same for both the controller and the observer. Thus, for the sake of brevity and to avoid repetitive material, we only show how to obtain the (full-state) control law, giving  the basic observer design and some additional remarks later in Section~\ref{sec-observer}.

\subsection{Design of a full-state feedback control law for unstable modes}
Based on the backstepping method \cite{krstic}, our idea is utilizing a invertible  Volterra integral transformation  
\begin{equation}\label{map}
	w_l^m(t,r)=u_l^m(t,r)-\int\limits\limits_{0}^{r} K_{lm}^{n}(r,\rho) u_l^m(t,\rho) d\rho,
\end{equation}
where the kernel $K_{lm}^{n}(r,\rho)$ is to be determined, which defined on the domain $\mathcal{T}_k = \{(r,\rho)\in\mathbb{R}^2 ;0 \leq \rho\leq r \leq R \}$ 
to convert the unstable system \eqref{eqn-un}-\eqref{eqn-un-bc} into an exponentially target system:
\begin{align}
	&  \partial_t w_l^m  =\epsilon \frac{\partial_r (r^{n-1} \partial_r w_l^m)}{r^{n-1}} - \epsilon l(l+n-2) \frac{w_l^m}{r^2} - c w_l^m ,\label{wl}\\
	& w_l^m(t,R) =0,\label{wl_bc}
\end{align}
where the constant $c>0$  is an adjustable convergence rate. From \eqref{map} and \eqref{wl_bc}, let $r=R$, we obtain the boundary control as the following full-state law
\begin{equation}\label{Un}
	U_l^m(t)=\int\limits\limits_{0}^{R} K_{lm}^{n}(R,\rho) u_l^m(t,\rho) d\rho.
\end{equation}
Following closely the steps of~\cite{nball} to find conditions for the kernels, and defining $K^n_{lm} (r,\rho)=G^n_{lm} (r,\rho)\rho\left(\frac{\rho}{r}\right)^{l+n-2}$, we finally reach a PDE that the $G$-kernels need to verify:
\begin{eqnarray}\label{eqn-gn}
\frac{\lambda(\rho)+c}{\epsilon} G^n_{lm} &=& \partial_{rr} G^n_{lm} +(3-n-2l)\frac{ \partial_r G^n_{lm}}{r}
-\partial_{\rho\rho} G^n_{lm}+(1-n-2l)\frac{ \partial_\rho G^n_{lm}}{\rho},\label{eqn-Gn}
\end{eqnarray}
with only one boundary condition:
\begin{eqnarray}
G^n_{lm}(r,r) 
&=& - \frac{\int_0^r (\lambda(\sigma)+c) d\sigma}{2r \epsilon }.\label{eqn-Gnbc}
\end{eqnarray}
We assume as usual that these  kernel equations are well-posed and the resulting kernel is bounded in $\mathcal{T}$; this will be analyzed later in Section~\ref{sec-wp}, providing also a numerical method for its computation.
\subsection{Closed-loop stability analysis of unstable modes}
To obtain the stability result of closed-loop system, we need three elements. We begin by stating the stability result for the target system. We follow by obtaining the existence of an inverse transformation that allows us to recover our original variable from the transformed variable. Then we relate the $L^2$ norm with spherical harmonics. With these elements, we construct the proof of stability mapping the result for the target system to the original system. This is done by showing that the transformation is an invertible map from $L^2$ into $L^2$.
 
We first discuss the stability of the target system, having the following Lemma:

\begin{lem} \label{Lem_wl}
	For all $l\in \mathbb{N}$, and for $c \geq 0$, the equilibrium  $w_l^m \equiv 0$ of system (\ref{wl})--(\ref{wl_bc}) is exponentially stable, i.e., there exists a positive constant $D_2$ such that for all $t$,
	\begin{eqnarray}
	\lVert w_l^m(t,\cdot) \lVert^2_{L^2} \leq e^{-D_2t} \lVert w_l^m(0,\cdot) \lVert^2_{L^2},
	\end{eqnarray}	
where the constant $D_2$ is independent of $n$, $l$ or $m$, and only depends on $c$, $\epsilon$, and $R$; it can be chosen as large as desired just by increasing the value of $c$.
\end{lem}
\begin{proof}
	Consider the Lyapunov function:
	\begin{equation}
	V_2(t)=\frac{1}{2}\lVert w_l^m(t,\cdot) \lVert_{L^2}^2,
	\end{equation} 
	then, taking its time derivative, we obtain
	\begin{align}
	\dot{V}_2  =& \int\limits_{0}^{R} \frac{\bar{w}_{l}^m \partial_t w_l^m + w_l^m \partial_t \bar{w}_{l}^m }{2} r^{n-1}dr
	\leq  - \left( \frac{\epsilon}{4R^2}  +c \right)\lVert w_l^m \lVert_{L^2}^2 
	\end{align}
	choosing 
	\begin{eqnarray}
	c=\frac{D_2}{2}-\frac{\epsilon}{4R^2}
	\end{eqnarray}
	we then obtain, independent of the value of $n$,
	\begin{eqnarray}
	\dot{V}_2 \leq -D_2 V_2,
	\end{eqnarray}
	thus proving the result.
\end{proof}

\begin{lem}\label{lem_ul_close}
	For $ \lvert l \lvert \leq L$,  let $c$ be chosen as in Lemma \ref{Lem_wl}, and assume that the kernel $K_{lm}^n(r,\rho)$ is bounded and integrable. The system \eqref{eqn-un} with boundary control \eqref{Un} is closed-loop exponentially stable, namely there exists positive constants $ C$ and $D_2$ such that
	\begin{equation}
	\lVert u_l^m(t,\cdot) \lVert^2_{L^2} \leq C e^{-D_2t} \lVert u_l^m(0,\cdot) \lVert^2_{L^2},
	\end{equation} 
	$C$ and $D_2$ are independent of $m$ or $l$, and only depend on $n$, $L$, $\lambda(r)$, $\epsilon$ and $R$.
\end{lem}
\begin{proof}
	The proof consists of two parts, one is existence of an inverse transformation, and then showing the equivalence of norms of the variables $u_{lm}^n$ and $w_{lm}^n$; the result then follows from the stability of the target system.

	As shown in~\cite{nball}, when $K_n(r,\rho)$ is bounded and integrable, the map \eqref{map} is reversible and its inverse transformation is
	\begin{eqnarray}\label{map_inverse_Kn}
	u_l^m(t,r)=w_l^m(t,r)+\int\limits_{0}^{r} L_{lm}^n(r,\rho) w_l^m(t,\rho) d\rho,
	\end{eqnarray}
	which is also bounded and integrable.
	Call now $\bar{K}$ and $\bar{L}$ the maximum of the bounds of the function $\check{K}_{lm}^{n}$ and $\check{L}_{lm}^{n}$ for a given $n$ and all $ l \leq L$ in their respective domains. It is easy to get
	\begin{eqnarray}
	&&\lVert w_l^m(t,\cdot) \lVert_{L^2}^2 \leq M_1 \lVert u_l^m (t,\cdot) \lVert_{L^2}^2, \label{ineq_w2u}\\
	&&\lVert u_l^m(t,\cdot) \lVert_{L^2}^2 \leq M_2 \lVert w_l^m(t,\cdot) \lVert_{L^2}^2.\label{ineq_u2w}
	\end{eqnarray}
	where $M_1=2+ R^4\bar{K}/(2n)$ and $M_2=2+ R^4\bar{L}/(2n)$.
	Combining  then Lemma \ref{Lem_wl} with the norm equivalence between $u_l^m$ and $w_l^m$ system stated as in \eqref{ineq_w2u} and \eqref{ineq_u2w}, it is easy to obtain
	\begin{align}
	&\lVert u_l^m(t,\cdot) \lVert^2_{L^2} \leq M_2 \lVert w_l^m(t,\cdot) \lVert^2_{L^2} 
	\leq M_2 e^{-D_2t} \lVert w_l^m(0,\cdot) \lVert^2_{L^2} \leq M_1 M_2 e^{-D_2t} \lVert u_l^m(0,\cdot) \lVert^2_{L^2}.
	\end{align}
	Let $C=M_1 M_2$, the result then follows.	
\end{proof}

Note that combining Lemmas~\ref{lem_ul_open} and~\ref{lem_ul_close} and taking $D=\min\{D_1,D_1\}$, we get the following stability result for \emph{all} spherical harmonics and thus the full physical system.

\begin{thm}\label{th-stab}
Under the assumption that  the kernel $K_{lm}^n(r,\rho)$ is bounded and integrable, the equilibrium $u_l^m\equiv 0$ of system \eqref{eqn-un}-\eqref{eqn-un-bc} under control law (\ref{Un})  is closed-loop exponentially stable, namely,  there exists a positive constant  $D$, such that for all $t$
	\begin{align}
		\label{equivalent1}\lVert u_l^m(t,\cdot) \lVert^2_{L^2} &\leq Ce^{-Dt}  \lVert u_l^m(0,\cdot) \lVert^2_{L^2}.
	\end{align}
where  $D$ can be chosen as large as desired just by increasing the value of $L$ and $c$ in the control design process.
\end{thm}

\section{Well-posedness of the kernel equations}\label{sec-wp}
Next, we state the main result of the paper, which was in part assumed in Theorem~\ref{th-stab}, also giving the requirements for $\lambda(r)$. In addition the proof of the result also provides a numerical method to compute the kernels, which is an alternative to successive approximations which do not work in this case (due to the singularities at the origin; see for instance~\cite{catalan} to see the resulting singular integral equation that needs to be solved).

\begin{thm}\label{th-Gresult}\rm
Under the assumption that $\lambda(r)$ is an even real analytic function in $[0,1]$, then for a given $n>1$ and all values of $l\in \mathbb{N}$, there is a unique power series solution $G^n_{lm}(r,\rho)$ for (\ref{eqn-gn})--(\ref{eqn-Gnbc}), even in its two variables in the domain $\mathcal{T}$, which is real analytic in the domain. In addition, if $\lambda(r)$ is analytic, but not even, then there is no power series solution to (\ref{eqn-gn})--(\ref{eqn-Gnbc}) for most values $l\in \mathbb{N}$.
\end{thm}

The requirement of evenness for $\lambda(r)$ might seem unusual. However, note that $r=\sqrt{x_1^2+x_2^2+\hdots+x_n^2}$, therefore in physical space $\lambda(r)$ will be non-smooth, unless it is even. Thus, while the kernels might exist for non-even $\lambda(r)$, we cannot expect them to be smooth, which might be indeed problematic for higher-order harmonics; not so much for lower-order, as shown in~\cite{catalan}, which only considers the 0-th order harmonic and consequently only requires boundedness of $\lambda(r)$.

\subsection{Proof of Theorem~\ref{th-Gresult}}
We start by giving out an algorithmic method to compute the power series for $G^n_{lm}(r,\rho)$, which will allow us to prove Theorem~\ref{th-Gresult} as well as numerically approximating the kernels.

First of all, we show that the evenness of $\lambda(r)$ is a necessary condition to find an analytic solution. Next, it is possible to establish that the series for $G^n_{lm}(r,\rho)$ only has even powers. Exploiting this property to suitably express  (\ref{eqn-gn})--(\ref{eqn-Gnbc}), we finally show the existence of the power series and thus Theorem~\ref{th-Gresult} follows. Convergence and related issues (radius of convergence) is studied towards the end, finishing the proof.

\subsubsection{Computing a power series solution for the kernels}\label{sec-series}
Starting from the most basic assumption of Theorem~\ref{th-Gresult}, we consider that $\lambda(r)$ is analytic in $[0,1]$, therefore it can be written as a convergent series (encompassing $c$ and $\epsilon$ for notational convenience):
\begin{equation}
\frac{\lambda(r)+c}{\epsilon}=\sum_{i=0}^\infty \lambda_i r^i, \label{def-powerlambda}
\end{equation}
which, by the evenness of $\lambda$, may only contain even powers\footnote{It is a known fact of analysis that even functions (respectively, odd functions) contain only even powers (respectively, odd powers) in their Taylor series. This fact has an easy proof by substituting the series in the definition of evenness $f(r)=f(-r)$ (respectively, oddness $f(r)=-f(-r)$) and checking the conditions verified by the coefficients.}, this is, $\lambda_i=0$ if $i$ is odd. We then seek for a solution of (\ref{eqn-gn})--(\ref{eqn-Gnbc}) of the form:
\begin{equation}
G^n_{lm}(r,\rho)=\sum_{i=0}^\infty \left(\sum_{j=0}^i C_{ij} r^j \rho^{i-j} \right), \label{def-power}
\end{equation}
where the dependence on $n$, $l$ and $m$ has been omitted for simplicity (the solution will depend on these values). The series in (\ref{def-power})  collects together (in the parenthesis) all the polynomial terms with the same degree.

It is easy to see that the boundary condition (\ref{eqn-Gnbc}) implies:
\begin{equation}
\forall i,\quad \sum_{j=0}^i C_{ij} =- \frac{\lambda_i}{2 (i+1)}, \label{eq-series1}
\end{equation}
which in particular implies $C_{00}=- \frac{\lambda_0}{2 \epsilon}$.
On the other hand, the left-hand side of (\ref{eqn-gn}) becomes
\begin{eqnarray}
\frac{\lambda(\rho)+c}{\epsilon} G^n_{lm} 
&=&\left[ \sum_{i=0}^\infty \left(\sum_{j=0}^i C_{ij} r^j \rho^{i-j} \right) \right] \sum_{i=0}^\infty \lambda_i \rho^i =\sum_{i=0}^\infty \left(\sum_{j=0}^i B_{ij} r^j \rho^{i-j} \right),
 \label{eq-seriesleft}
\end{eqnarray}
where we have defined 
\begin{equation}
B_{ij}=\sum_{k=j}^i C_{kj} \lambda_{i-k}. \label{def-Bij}
\end{equation}
Finally, to express the right-hand side of (\ref{eqn-gn}), denote $\gamma=n+2l-2\geq 0$ and define the operators $D_1=\partial_{rr} +(1-\gamma)1/r \partial_r$ and $D_2=-\partial_{\rho\rho} +(-1-\gamma)1/\rho  \partial_\rho$. Then
\begin{eqnarray}
D_1 G^n_{lm}&=&\sum_{i=1}^\infty \left(\sum_{j=1}^i j(j-\gamma)C_{ij} r^{j-2} \rho^{i-j} \right), \label{eqn-D1} \\
D_2 G^n_{lm}&=&\sum_{i=1}^\infty \left(\sum_{j=0}^{i-1} (i-j)(j-i-\gamma)C_{ij} r^j \rho^{i-j-2} \right)\hspace{-2pt},\quad ~~~ \label{eqn-D2} 
\end{eqnarray}
and thus, rewriting the sum to be homogeneous with (\ref{eq-seriesleft}),we find
\begin{equation}
(D_1+D_2)G^n_{lm}=\sum_{i=-1}^\infty \left(\sum_{j=-1}^{i+1} D_{ij} r^j \rho^{i-j} \right),\label{eqn-DD}
\end{equation}
where, (assuming $C_{ij}=0$ if $i,j<0$ or $j>i$),
\begin{eqnarray}
D_{ij}&=&  (j+2)(j+2-\gamma)C_{(i+2)(j+2)} -(i-j+2)(i-j+2+\gamma)C_{(i+2)j},
\end{eqnarray}

Equating (\ref{eqn-DD}) and (\ref{eq-seriesleft}), we obtain a system of equations:
\begin{eqnarray}
\forall i\geq -1,&& \quad D_{i(i+1)}=D_{i(-1)}=0, \\
\forall i\geq 0, 0\leq j \leq i,&& \quad  (j+2)(j+2-\gamma)C_{(i+2)(j+2)}-(i-j+2)(i-j+2+\gamma)C_{(i+2)j} =B_{ij}.\quad
\end{eqnarray}

With $\lambda(r)$ and $n$ are fixed, we want to show that the kernel equations are solvable for all values of $l\in\mathbb{N}$. Thus, $\gamma$ takes  increasing values. In addition we can assume $\gamma \neq 1$, since the case $n=3$, $l=0$ was already addressed in \cite{sphere} showing that it reduces to the usual 1-D kernel equations for parabolic systems~\cite{krstic}, which admits a power series solution according to~\cite{ascencio}.

The first two equalities, if $\gamma \neq 1$,  imply
\begin{equation}
\forall i\geq 1, \quad C_{i1}=C_{i(i-1)}=0,\label{eq-C1s}
\end{equation}
and, in particular, $C_{10}=C_{01}=0$, whereas the second equality results in a system of equations that needs to be solved recursively, starting at $i=0$. It can be rewritten as follows to start at $i=2$ (since $C_{00}$, $C_{10}$ and $C_{01}$ are already determined).
\begin{equation}
\forall i\geq2,\,0\leq j \leq i-2, \,  (j+2)(j+2-\gamma)C_{i(j+2)} -(i-j)(j-i-\gamma)C_{ij} =
\sum_{k=j}^{i-2} C_{kj} \lambda_{i-2-k}. \label{eq-recursive}
\end{equation}
Note that for each $i\geq 2$, there are $i+1$ coefficients in (\ref{def-power}) but $i+2$ relations: one from (\ref{eq-series1}), two from (\ref{eq-C1s}) and $i-1$ from (\ref{eq-recursive}).  Thus, it would seem  that (\ref{eq-series1})--(\ref{eq-C1s})--(\ref{eq-recursive}) is in general an incompatible system. This is indeed the case if $\lambda(r)$ is not even, i.e., if the series (\ref{def-powerlambda}) contains odd powers, as shown in the next section.
\subsubsection{Evennes requirement of $\lambda(r)$}
We start with the following result.
\begin{lem}\label{lem1}\rm
If $\lambda(r)$ is not even, then there are values of $l\in\mathbb{N}$ for which there is no  solution to (\ref{eqn-gn})--(\ref{eqn-Gnbc}) in the form of (\ref{def-power}).
\end{lem}
\begin{proof}
We show that, if there exists $i$ odd such that $\lambda_i\neq0$, then there is no solution in the form of a power series. First, if $\lambda_1 \neq 0$, then from (\ref{eq-series1}) we know that $C_{01}+C_{10}=- \frac{\lambda_1}{4\epsilon}$, however since form (\ref{eq-C1s}) one has $C_{01}=C_{10}=0$, this cannot hold. Consider now there is indeed a value $i>1$ for which a coefficient $\lambda_{i}$  is distinct from zero and let us show the result by contradiction. Consider the first such $i$. Now, since the right-hand side of (\ref{eq-recursive}) depends on $C_{(i-2)j}$, one gets that for all odd $i'<i$ $C_{i'j}$ must zero from  (\ref{eq-series1})--(\ref{eq-C1s})--(\ref{eq-recursive}) all having a zero right-hand side (this can be formalized with an induction argument; we skip the details). Thus, at $i$, the following system of equations has to be verified:
\begin{eqnarray}
C_{i1}&=&C_{i(i-1)}=0,\label{eq-start} \\
\sum_{j=0}^i C_{ij} &=&- \frac{\lambda_i}{2 \epsilon(i+1)}, \label{eq-sumrequired}
\end{eqnarray}
and for $0\leq j \leq i-2$,
\begin{equation}
(j+2)(j+2-\gamma)C_{i(j+2)} -(i-j)(i-j+\gamma)C_{ij}=0,\label{eq-sequenceequalities}
\end{equation}
Let us consider $l$ sufficiently large such that $\gamma>i$, so that the coefficient $(j+2-\gamma)$ in (\ref{eq-sequenceequalities}) is distinct from zero in the full range of $j$, namely $0\leq j \leq i-2$. Then none of the coefficients in (\ref{eq-sequenceequalities}) is zero.
Therefore, combining (\ref{eq-start}) with (\ref{eq-sequenceequalities}), from $C_{i1}$ we can find $C_{i3}$, then $C_{i5}$, and so on. Similarly, from $C_{i(i-1)}$ we can find $C_{i(i-3)}$, $C_{i(i-5)}$ and so on. These two sequences don't overlap because $i$ is odd and therefore, one finds $C_{ij}=0$ for all $0\leq j \leq i$ which is not compatible with (\ref{eq-sumrequired}) unless $\lambda_{i}=0$, which contradicts our initial assumption.
\end{proof}
Next we show that evenness of $\lambda$ implies evenness of the kernels.
\begin{lem}\label{lem2}\rm
If $\lambda(r)$ is even, then, a solution to (\ref{eqn-gn})--(\ref{eqn-Gnbc}) in the form of (\ref{def-power})  only has even powers.
\end{lem}
\begin{proof}
We need to prove that $C_{ij}=0$ if either $i$ or $j$ is odd.
From the proof of Lemma~\ref{lem1}, we directly know that for odd $i$ one has $C_{ij}=0$.  Fix, then, $i$ even and consider $j$ odd;  for $i=2$, the result is obvious. Assuming $C_{i'j}=0$ for all even numbers $i'<i$ and $j$ odd, let us prove the result by induction on the first coefficient. As before, we would need to solve (\ref{eq-sumrequired})--(\ref{eq-recursive}). The right-hand side $B_{ij}=\sum_{k=j}^{i-2} C_{kj} \lambda_{i-2-k}$ of (\ref{eq-recursive}) is zero as in (\ref{eq-sequenceequalities}) by the induction hypothesis (if $k$ even) or directly zero if $k$ odd. Then, following again the proof of Lemma~\ref{lem1}, we have the same system of equations (\ref{eq-sumrequired})--(\ref{eq-sequenceequalities}) for our even $i$ and odd $j$'s. Now:
$$C_{ij}=\frac{(j+2)(j+2-\gamma)}{(i-j)(i-j+\gamma)}C_{i(j+2)},$$
so starting from $C_{i(i-1)}=0$ we find $C_{i(i-3)}=0$, then $C_{i(i-5)}$, and so on; however, with $i$ being even, this sequence ends now in $C_{i1}$ (thus, the proof of Lemma~\ref{lem1} does not apply because the sequences starting at $C_{i1}$ and $C_{i(i-1)}$ overlap). Thus, one finds $C_{ij}=0$ for all odd values of $j$ between 1 and $i-1$. 
\end{proof}
\subsubsection{Well-posedness of the coefficient system}
Next, we show that the coefficients of the power series can always be found, which by the previous Lemmas only requires studying the even coefficients. For simplification, we redefine (\ref{def-powerlambda}) and (\ref{def-power}) as:
\begin{equation}
\frac{\lambda(r)+c}{\epsilon}=\sum_{i=0}^\infty \lambda_i r^{2i}, \quad
G^n_{lm}(r,\rho)=\sum_{i=0}^\infty \left(\sum_{j=0}^i C_{ij} r^{2j} \rho^{2(i-j)} \right), \label{def-power2}
\end{equation}
without bothering to redefine the coefficients (note that (\ref{def-Bij}) does not require any change). Defining as well $\gamma'=\frac{\gamma}{2}=\frac{n}{2}+l-1\geq 0$, the new system of equations to be solved is
\begin{eqnarray}
\forall i,\quad && \sum_{j=0}^i C_{ij} =- \frac{\lambda_i}{2 (2i+1)}, \label{eq-sumrequired1} 
\end{eqnarray}
and
\begin{equation}
\forall i\geq1,0\leq j \leq i-1, \,\,  (j+1)(j+1-\gamma')C_{i(j+1)} +(i-j)(j-i-\gamma')C_{ij} =
\sum_{k=j}^{i-1} C_{kj} \lambda_{i-1-k}=B_{{(i-1})j}. \quad\label{eq-recursive1}
\end{equation}
Let us outline the solution procedure, and later derive some conclusions.  Solving in (\ref{eq-recursive1}) every $C_{ij}$ as a function of $C_{i(j+1)}$ we  get:
\begin{equation}
 C_{ij}=\frac{(j+1)(j+1-\gamma') C_{i(j+1)}+ B_{(i-1)j}}{(i-j)(i-j+\gamma')}, \label{eq-Cij}
\end{equation}
which can be written more briefly if we define, for $i>0$ and $0\leq j<i$,
\begin{equation}
a_{ij}(\gamma')=\frac{(j+1)(j+1-\gamma')}{(i-j)(i-j+\gamma')},\label{eq-aij}
\end{equation}
as
\begin{equation}
C_{ij}=a_{ij}(\gamma')C_{i(j+1)}+\frac{B_{(i-1)j}}{(i-j)(i-j+\gamma')}.
\end{equation}
To be able to simplify a bit the equation, redefine
\begin{equation}
\hat B_{(i-1)j}=\frac{B_{(i-1)j}}{(i-j)(i-j+\gamma')}
\end{equation}
then,
\begin{equation}
C_{ij}=a_{ij}(\gamma')C_{i(j+1)}+\hat B_{(i-1)j}.
\end{equation}
and iterating this equality until reaching $C_{ii}$, we get
\begin{eqnarray*}\label{eq-coefs}
 C_{ij}&=&\left[\prod_{k=j}^{k=i-1} a_{ik}(\gamma')  \right]C_{ii}  +  \hat B_{(i-1)j} + \sum_{r=j+1}^{i-1}\prod_{k=r}^{k=i-1} a_{ik}(\gamma') \hat B_{(i-1)r} ,
\end{eqnarray*}
and inserting this into (\ref{eq-sumrequired1}), we reach an equation for $ C_{ii}$, namely
\begin{equation}
 C_{ii}=- \frac{1}{\kappa(i,\gamma')}\left[ \frac{\lambda_{i}}{2 \epsilon(2i+1)}+H_{i}\right], \label{eq-Cii}
\end{equation}
where
\begin{equation}
\kappa(i,\gamma')=1+\sum_{j=0}^{i-1} \prod_{k=j}^{k=i-1} a_{ik}(\gamma'), \label{eq-kappa-aij}
\end{equation}
and
\begin{eqnarray}
H_i&=&\sum_{j=0}^{i-1} 
 \hat B_{(i-1)j} +\sum_{j=0}^{i-2} 
\sum_{r=j+1}^{i-1} \prod_{k=r}^{k=i-1} a_{ik}(\gamma') \hat B_{(i-1)r} 
 \nonumber \\ 
&=&\sum_{j=0}^{i-1} 
 \hat B_{(i-1)j} +\sum_{r=1}^{i-1}  \sum_{j=0}^{r-1}  \prod_{k=r}^{k=i-1} a_{ik}(\gamma') \hat B_{(i-1)r} \nonumber \\ 
&=&\sum_{j=0}^{i-1} 
 \hat B_{(i-1)j}+\sum_{r=1}^{i-1}  r  \prod_{k=r}^{k=i-1} a_{ik}(\gamma') \hat B_{(i-1)r}
 \nonumber \\ 
&=&\sum_{j=0}^{i-1} \left(1+j  \prod_{k=j}^{k=i-1} a_{ik}(\gamma') \right) \hat B_{(i-1)j}
\end{eqnarray}
It is quite clear that these $\kappa(i,\gamma')$ will play an important role; in particular, if they are non-zero, one can always find a unique solution for the coefficients $C_{ij}$. Thus one needs to show that $\kappa(i,\gamma')\neq0$ for any possible $i$ or $\gamma'$. The following lemma shows this is indeed the case, by exploiting a connection of the $a_
{ij}$ coefficients with Gauss' hypergeometric functions.
%
%
\begin{lem}\rm\label{lem-gauss}
Let $i$ be a positive integer and $\gamma'\geq 0$ a real number. Then, it holds that 
%
\begin{eqnarray}
\kappa(i,\gamma')=
\frac{2i!}{i!}\frac{\Gamma({\gamma'+1})}{\Gamma({i+\gamma'+1})}>0,\label{eq-kappa-formula-exp3}\end{eqnarray}
where $\Gamma$ denotes the Gamma function~\cite[p.255]{abramowitz}.
\end{lem}
\begin{proof}
Recalling from (\ref{eq-kappa-aij}) and (\ref{eq-aij}) the definitions of $\kappa(i,\gamma')$ and $a_{ij}$, respectively, one has
\begin{equation}
\kappa(i,\gamma')=1+\sum_{j=0}^{i-1} \prod_{k=j}^{k=i-1} \frac{(k+1)(k+1-\gamma')}{(i-k)(i-k+\gamma')}
\end{equation}
which can be rewritten in terms of binomial numbers and rising/falling factorials\footnote{Rising factorials $(x)^{\overline n}$ are sometimes expressed using the Pochhammer's symbol, with a slightly different notation, namely $(x)_n$.}~\cite{knuth} as
\begin{equation}
\kappa(i,\gamma')=\sum_{j=0}^{i}
\left(\begin{array}{c} i \\ j\end{array}\right)
\frac{(i-\gamma')^{\underline{i-j}}}{(1+\gamma')^{\overline{i-j}}}
\end{equation}
and reordering the sum and using $\left(\begin{array}{c} i \\ j\end{array}\right)=\left(\begin{array}{c} i \\ i-j\end{array}\right)$,
\begin{equation}
\kappa(i,\gamma')=\sum_{j=0}^{i}
\left(\begin{array}{c} i \\ j\end{array}\right)
\frac{(i-\gamma')^{\underline j}}{(1+\gamma')^{\overline j}}= 
\sum_{j=0}^{i} (-1)^j
\left(\begin{array}{c} i \\ j\end{array}\right)
\frac{(\gamma'-i)^{\overline j}}{(1+\gamma')^{\overline j}},
\end{equation}
where the obvious fact $(x)^{\underline j}=(-1)^j (-x)^{\overline j}$ has been used. 
Consider now the finite polynomial $p_i(x,\gamma')$ defined as
\begin{equation}
p_i(x,\gamma)=
\sum_{j=0}^{i} (-1)^j
\left(\begin{array}{c} i \\ j\end{array}\right)
\frac{(\gamma'-i)^{\overline j}}{(1+\gamma')^{\overline j}} x^j
\end{equation}
From the definition of Gauss' hypergeometric function~\cite[p.561]{abramowitz}, denoted as ${}_2 F_1(a,b;c;x)$, in the polynomial case ($a$ or $b$ non-positive integer) and noting $(-1)^j
\left(\begin{array}{c} i \\ j\end{array}\right)=\frac{(-i)^{\overline j}}{j!}$, it is immediate that 
\begin{equation}
p_i(x,\gamma)={}_2 F_1(-i,\gamma'-i;1+\gamma';x)
\end{equation}
and therefore, from Gauss' summation theorem~\cite[p.556]{abramowitz}, which is applicable in this case since $1+2i>0$,
\begin{equation}
\kappa(i,\gamma')=p_i(1,\gamma)={}_2 F_1(-i,\gamma'-i;1+\gamma';1)=\frac{\Gamma(1+\gamma')}{\Gamma(1+\gamma'+i)}\frac{\Gamma(1+2i)}{\Gamma(1+i)}=\frac{2i!}{i!} \frac{\Gamma(1+\gamma')}{\Gamma(1+\gamma'+i)},
\end{equation}
finishing the proof.
\end{proof}

The next result is an immediate conclusion of the positivity of $\kappa(i,\gamma')$:
\begin{lem}\rm
If $\lambda(r)$ is even, then, for all values of $l\in\mathbb{R}$, the coefficientes in (\ref{def-power2}) that solve (\ref{eqn-gn})--(\ref{eqn-Gnbc}) can be uniquely found up to any order $i$.
\end{lem}

To conclude the proof of Theorem~\ref{th-Gresult}, we need to prove analyticity of the series (\ref{def-power2}). This step, however, requires splitting the problem in two possible cases: odd dimension (thus, $\gamma'=n/2+2l-1$ is not an integer) and even dimension ($\gamma'$  integer).
\subsubsection{Proof of analyticity for odd dimension}\label{sec-odd}
In the odd-dimension case, define the following coefficients:
\begin{equation}
L_{i0}=1, \qquad 
L_{ij}=\left(\begin{array}{c} i \\ j\end{array}\right) \frac{(i+\gamma')(i-1+\gamma')\hdots(i-j+\gamma'+1)}{(1-\gamma')(2-\gamma')\hdots(j-\gamma')}, \quad j>0
\label{def-Lij}
\end{equation}
with $L_{i0}$ defined as 1; these are well-defined given that $\gamma'$ is non-integer. They can also be expressed as
$$
L_{ij}=\left(\begin{array}{c} i \\ j\end{array}\right) \frac{\Gamma(1-\gamma')\Gamma(i+1+\gamma')}{\Gamma(j+1-\gamma')\Gamma(i-j+1+\gamma')}
$$
Now, in (\ref{eq-sumrequired1})--(\ref{eq-recursive1}), denote $ C_{ij}=L_{ij}\check C_{ij}$. Replacing in the recurrence we get
\begin{eqnarray}
B_{{(i-1})j}&=&(j+1)(j+1-\gamma')
\left(\begin{array}{c} i \\ j+1 \end{array}\right) \frac{\Gamma(1-\gamma')\Gamma(i+1+\gamma')}{\Gamma(j+2-\gamma')\Gamma(i-j+\gamma')}\check C_{i(j+1)} 
\nonumber \\&&
-(i-j)(i-j+\gamma')\left(\begin{array}{c} i \\ j\end{array}\right) \frac{\Gamma(1-\gamma')\Gamma(i+1+\gamma')}{\Gamma(j+1-\gamma')\Gamma(i-j+1+\gamma')}
\check C_{ij} 
\nonumber \\
&=&
(i-j)(i-j+\gamma')L_{ij}
\left(\check C_{i(j+1)} -\check C_{ij} \right)
\end{eqnarray}
Define now
\begin{equation}
\check 
B_{{(i-1})j}
=\frac{
B_{{(i-1})j}}{ 
(i-j)(i-j+\gamma')L_{ij}
}\label{eqn-checkB}
\end{equation}
and the new set of recurrence equations for $\check C_{ij}$ becomes rather simple:
\begin{eqnarray}
\forall i,\quad && \sum_{j=0}^i L_{ij} \check C_{ij} =- \frac{\lambda_i}{2 (2i+1)},\label{eqn-sumrequiredcheck} \\
\forall i\geq1,0\leq j \leq i-1,&& \check C_{i(j+1)} - \check C_{ij} =
\check B_{{(i-1})j},
\end{eqnarray}
and the recurrence is easily solvable in terms of one element; for instance, $\check C_{ii}$:
\begin{equation}
\check C_{ij}=\check C_{ii} - \sum_{r=j}^{r=i-1} \check B_{{(i-1})r} \label{eq-recursioncii}
\end{equation}
for $j=1,\hdots,i$.
Replacing in (\ref{eqn-sumrequiredcheck})  we reach
$$
 \sum_{j=0}^i L_{ij} \check C_{ii} -
  \sum_{j=0}^{i-1} L_{ij} \sum_{r=j}^{r=i-1} \check B_{{(i-1})r}
  =- \frac{\lambda_i}{2 (2i+1)}
$$
Thus:
\begin{eqnarray}
\check C_{ii} &=& \frac{ -\frac{\lambda_i}{2 (2i+1)} +  \sum_{k=0}^{k=i-1} \sum_{r=k}^{r=i-1}L_{ik}  \check B_{{(i-1})r} }{ \sum_{k=0}^{k=i} L_{ik} }\nonumber \\
&=&\frac{- \frac{\lambda_i}{2 (2i+1)} +  \sum_{r=0}^{r=i-1} \check B_{{(i-1})r} \left( \sum_{k=0}^{k=r}L_{ik} \right) }{  \sum_{k=0}^{k=i} L_{ik}  }
\end{eqnarray}
Call 
\begin{equation}
R_{ij}= \sum_{k=0}^{k=j}L_{ik},\qquad R_i=\sum_{k=0}^{k=i}L_{ik}=R_{ii}.\label{def-Rij}
\end{equation}
Then
\begin{eqnarray}
\check C_{ii} 
&=&\frac{ -\frac{\lambda_i}{2 (2i+1)} +  \sum_{r=0}^{r=i-1} \check B_{{(i-1})r}R_{ir} }{  R_i  }
\end{eqnarray}
Now, solving for the remaining coefficients from (\ref{eq-recursioncii}):
\begin{eqnarray}
\check C_{ij}&=&\frac{ -\frac{\lambda_i}{2 (2i+1)} +  \sum_{r=0}^{r=i-1} \check B_{{(i-1})r}R_{ir} }{  R_i  } - \sum_{r=j}^{r=i-1} \check B_{{(i-1})r} \nonumber\\ \nonumber  &=&\frac{ -\frac{\lambda_i}{2 (2i+1)} +  \sum_{r=0}^{r=i-1} \check B_{{(i-1})r}R_{ir} - \sum_{r=j}^{r=i-1} \check B_{{(i-1})r} R_i}{  R_i  }  \nonumber \\
&=&
\frac{- \frac{\lambda_i}{2 (2i+1)} +  \sum_{r=0}^{r=j-1} \check B_{{(i-1})r}R_{ir} - \sum_{r=j}^{r=i-1} \check B_{{(i-1})r} (R_i-R_{ir})}{  R_i  }  
\end{eqnarray}
Finally, recovering the coefficients $C_{ij}$ from $C_{ij}=L_{ij}\check C_{ij}$ and using (\ref{eqn-checkB}):
\begin{eqnarray}
 C_{ij}&=&
\frac{ -\frac{\lambda_i}{2 (2i+1)}L_{ij} +  \sum_{r=0}^{r=j-1} B_{{(i-1})r} \frac{L_{ij} R_{ir}
}{ 
(i-r)(i-r+\gamma')L_{ir}} - \sum_{r=j}^{r=i-1} B_{{(i-1})r} \frac{
L_{ij}(R_i-R_{ir})}{ 
(i-r)(i-r+\gamma')L_{ir}}}{  R_i  },
\end{eqnarray}
which is quite explicit.

Notice that since $\rho \leq r$,
$$\vert G^n_{lm}(r,\rho) \vert \leq \sum_{i=0}^\infty \left(\sum_{j=0}^i \vert C_{ij} \vert r^{2j} \rho^{2(i-j)} \right) 
\leq \sum_{i=0}^\infty r^{2i} \left(\sum_{j=0}^i \vert C_{ij}  \right), $$
thus, defining $\alpha_i=\sum_{j=0}^i \vert C_{ij} \vert$, if we can prove that  $ \sum_{i=0}^\infty \alpha_i r^{2i} $ converges for a certain radius of convergence $R$, so does $G^n_{lm}(r,\rho)$ for $\rho\leq r \leq R$, and thus we obtain the required analyticity.
Now:
\begin{eqnarray}
\alpha_i&=& \sum_{j=0}^i \vert C_{ij}\vert \nonumber \\ &\leq &
 \frac{\vert \lambda_i \vert }{2 (2i+1)} \frac{ \sum_{j=0}^i \vert L_{ij} \vert}{\vert R_i \vert} +\frac{  \sum_{j=0}^i\sum_{r=0}^{r=j-1} \vert B_{{(i-1})r} \vert\vert \frac{L_{ij} R_{ir}
}{ 
(i-r)(i-r+\gamma')L_{ir}}\vert }{  \vert R_i \vert }
 \nonumber \\ &&
+\frac{\sum_{j=0}^i \sum_{r=j}^{r=i-1} \vert B_{{(i-1})r} \vert\vert\frac{
L_{ij}(R_i-R_{ir})}{ 
(i-r)(i-r+\gamma')L_{ir}}\vert}{  \vert R_i \vert }
\nonumber \\ &= &
 \frac{\vert \lambda_i \vert }{2 (2i+1)} \frac{ \sum_{j=0}^i \vert L_{ij} \vert}{\vert R_i \vert} +\frac{  \sum_{r=0}^{i-1}\sum_{j=r+1}^{i} \vert B_{{(i-1})r} \vert\vert \frac{L_{ij} R_{ir}
}{ 
(i-r)(i-r+\gamma')L_{ir}}\vert }{  \vert R_i \vert }
 \nonumber \\ &&
+\frac{\sum_{r=0}^{r=i-1}\sum_{j=0}^r \vert B_{{(i-1})r} \vert\vert\frac{
L_{ij}(R_i-R_{ir})}{ 
(i-r)(i-r+\gamma')L_{ir}}\vert}{  \vert R_i \vert }
\nonumber \\ &\leq &
 \frac{\vert \lambda_i \vert }{2 (2i+1)} \frac{ \sum_{j=0}^i \vert L_{ij} \vert}{\vert R_i \vert} +\sum_{r=0}^{i-1}\vert B_{{(i-1})r} \vert \left(  \frac{  \vert R_{ir} \vert \sum_{j=r+1}^{i} \vert L_{ij} \vert + \vert R_i-R_{ir}\vert \sum_{j=0}^r  \vert L_{ij} \vert}{ 
(i-r)(i-r+\gamma')\vert L_{ir} \vert \vert R_i \vert }\right) \label{eqn-alphandef}
,
\end{eqnarray}
To prove the convergence of the power series $\sum_{i=0}^\infty \alpha_i r^{2i}$ consider the following lemma.
\begin{lem}\label{lem-conv}
Consider $g(x)=\sum_{i=0}^\infty g_i x^{2i}$ and $h(x)=\sum_{i=0}^\infty h_i x^{2i}$ analytic functions, both with radius of convergence $R$.
Let $i_0$ be a nonnegative integer, let $(a_i)_{i=0}^\infty$ be a sequence of real numbers, and define $f(x)=\sum_{i=0}^\infty a_i x^{2i}$, where $a_i$ verify, for $i>i_0$
\begin{equation}
a_i \leq  b_i \vert g_i \vert + c_i \sum_{j=0}^{i-1} a_j \vert h_{i-1-j} \vert
\end{equation}
where the sequences $b_i,c_i \geq 0$ are decreasing for $i>i_0$, with $c_i$ also verifying $\lim_{i \rightarrow \infty} c_i=0$.
Then, $f(x)$ is analytic with radius of convergence at least $R$.
\end{lem}
\begin{proof}
Since $g$ and $h$ analytic with radius of convergence $R$ we can write $\vert g_i\vert,\vert h_i \vert \leq  M R^{-2i}$, where the definition as  power series of squares has been taken into account.
Thus:
\begin{equation}
a_i  \leq b_i M R^{-2i} + c_i \sum_{j=0}^{i-1} a_j M R^{-2i+1+2j}
\end{equation}
Define $\check a_i=a_i$ for $i\leq i_0$ and, for $i> i_0$, $\check a_i  = b_i M R^{-2i} + c_i \sum_{j=0}^{i-1} \check a_j M R^{-2i+1+2j}$. Obvioulsy $a_i \leq \check a_i$ and therefore the radius of convergence of $f(x)$ would be at least the radius of convergence of $\check f(x)=\sum_{i=0}^\infty \check a_i x^{2i}$. 
Now:
\begin{equation}
\check a_{i+1}  = b_{i+1} M R^{-2i-2} + c_{i+1} \sum_{j=0}^{i} \check a_j M R^{-2i+2j-1} = b_{i+1} M R^{-2i-2} + c_{i+1} M \check a_i+c_{i+1} R^{-2}\sum_{j=0}^{i-1} \check a_j M R^{-2i2+j+1} 
\end{equation}
It is sufficient to compute
\begin{eqnarray}
\lim_{i \rightarrow \infty} \frac{\check a_{i+1}}{\check a_i} 
 &=& \lim_{i \rightarrow \infty} \frac{ b_{i+1} M R^{-2i-2} + c_{i+1} M \check a_i+ c_{i+1}R^{-2}\sum_{j=0}^{i-1} \check a_j M R^{-2i+2j+1} }{\check a_i}
\nonumber \\ &=& \lim_{i \rightarrow \infty} 
M c_{i+1} +  \lim_{i \rightarrow \infty}  \frac{ b_{i+1} M R^{-2i-2} +c_{i+1} R^{-2}\sum_{j=0}^{i-1} \check a_j M R^{-2i+2j+1} }{b_i M R^{-2i} + c_i \sum_{j=0}^{i-1} \check a_j M R^{-2i+1+2j}}
\nonumber \\ &= &R^{-2} \lim_{i \rightarrow \infty}  \frac{ b_{i+1} +c_{i+1} \sum_{j=0}^{i-1} \check a_j R^{2j+1} }{b_i   + c_i \sum_{j=0}^{i-1} \check a_j R^{1+2j}}
\nonumber \\ &\leq &R^{-2} \lim_{i \rightarrow \infty}  \frac{ b_{i}+c_i \sum_{j=0}^{i-1} \check a_j R^{2j+1} }{b_i   + c_i \sum_{j=0}^{i-1} \check a_j R^{1+2j}}
\nonumber \\ &=&R^{-2},
\end{eqnarray}
where the inequality holds for sufficiently large $i>i_0$ and thus in the limit, therefore proving the lemma.
\end{proof}
To apply Lemma~\ref{lem-conv} to (\ref{eqn-alphandef}) we need to bound some of the terms. In particular, if we are able to find $b_i$ and $c_i$ such that \begin{equation}
\frac{ \sum_{j=0}^i \vert L_{ij} \vert}{2 (2i+1) \vert R_i \vert} \leq b_i, \quad 
\max_{r\in\{0,\hdots,i-1\}} \left(  \frac{  \vert R_{ir} \vert \sum_{j=r+1}^{i} \vert L_{ij} \vert + \vert R_i-R_{ir}\vert \sum_{j=0}^r  \vert L_{ij} \vert}{ 
(i-r)(i-r+\gamma')\vert L_{ir} \vert \vert R_i \vert }\right) \leq c_i, \label{eq-bounds}
\end{equation}
we get
\begin{eqnarray}
\alpha_i &\leq & b_i \vert \lambda_i \vert  + c_i \sum_{r=0}^{i-1}\vert B_{{(i-1})r} \vert 
\nonumber \\ &=& 
 b_i \vert \lambda_i \vert  + c_i \sum_{r=0}^{i-1}
\sum_{k=r}^{i-1} \vert C_{kr} \vert \vert \lambda_{i-1-k}  \vert
\nonumber \\ &=& 
 b_i \vert \lambda_i \vert  + c_i \sum_{k=0}^{i-1}
\sum_{r=0}^{k} \vert C_{kr} \vert \vert \lambda_{i-1-k}  \vert
\nonumber \\ &=& 
 b_i \vert \lambda_i \vert  + c_i \sum_{k=0}^{i-1}
\alpha_k \vert \lambda_{i-1-k}  \vert
,
\end{eqnarray}
so, assuming $b_i$ and $c_i$ verify the conditions given in Lemma~\ref{lem-conv}, and given that $\lambda(x)$ has a radius of convergence of at least one, we obtain that $G^n_{lm}(r,\rho)$ converges and defines an analytic function for $\rho\leq r \leq 1$, thus proving Theorem~\ref{th-Gresult} for the odd dimension case.

It remains to find such $b_i$ and $c_i$. Proceeding exactly as in Lemma~\ref{lem-gauss} with a slight modification, we directly find
\begin{equation}R_i={}_2 F_1(-i,-(\gamma'+i);1-\gamma';1)=\frac{\Gamma(1-\gamma')\Gamma(2i+1)}{\Gamma(-\gamma'+1+i)\Gamma(i+1)}=\frac{2i!}{i!} \frac{\Gamma(1-\gamma')}{\Gamma(1-\gamma'+i)} \label{eqn-explicitRi}
\end{equation}


Now, let $N=\gamma'-1/2$ and $i>2N$. One can see that for $1\leq j \leq N$,
\begin{eqnarray}
\vert L_{ij} \vert &=& \left(\begin{array}{c} i \\ j\end{array}\right) \frac{(i+\gamma')(i-1+\gamma')\hdots(i-j+\gamma'+1)}{(\gamma'-1)(\gamma'-2)\hdots(\gamma'-j)} \nonumber \\
&\leq&\vert L_{iN} \vert
\nonumber \\
&=&
 \left(\begin{array}{c} i \\ N \end{array}\right) \frac{(i+\gamma')(i-1+\gamma')\hdots(i+3/2)}{(\gamma'-1)(\gamma'-2)\hdots(1/2)}
\nonumber \\
& \leq &\frac{i!}{N! (i-N)!} \frac{(i+N+1)!}{(i+1)!}\frac{1}{(N-1)!}
\end{eqnarray}
and
\begin{equation}
\vert R_i\vert =\frac{2i!}{i!} \left| \frac{ \Gamma(\gamma'-1)  }{\Gamma(i+1-\gamma')} \right| \geq \frac{2i!}{i! (i-N)!N!}. \label{eqn-Ribound}
\end{equation}
Thus, for $i>2N$, calling $d_i$ the following sequence
\begin{equation}
d_i=\frac{ \sum_{j=0}^N \vert L_{ij} \vert }{ \vert R_i \vert} \leq
\frac{i! (i+N+1)!}{(i+1) (N-1)! 2i!}  \label{eqn-di}
\end{equation}
it is clear that $d_i$ is a decreasing sequence, since from the ratio test  $r=\lim_{i\rightarrow \infty} \frac{d_{i+1}}{d_i}=1/4$.

Now, set $i_0=2N$. For $i>i_0$,
\begin{equation}
\frac{ \sum_{j=0}^i \vert L_{ij} \vert}{2 \epsilon(2i+1) \vert R_i \vert} 
=\frac{ \sum_{j=0}^N \vert L_{ij} \vert+\sum_{j=N+1}^i \vert L_{ij} \vert }{2 (2i+1) \vert R_i \vert} \leq \frac{ d_i +1}{2 (2i+1)} =b_i
\end{equation}
It is obvious that $b_i$ is decreasing since $d_i$ is decreasing.

Now we need to find a sequence $c_i$ for the second term in (\ref{eq-bounds}). First of all,
\begin{equation}
  \frac{  \vert R_{ir} \vert \sum_{j=r+1}^{i} \vert L_{ij} \vert + \vert R_i-R_{ir}\vert \sum_{j=0}^r  \vert L_{ij} \vert}{ 
(i-r)(i-r+\gamma')\vert L_{ir} \vert \vert R_i \vert } \leq   \frac{  2   \sum_{j=r+1}^{i} \vert L_{ij} \vert   \sum_{j=0}^r  \vert L_{ij} \vert}{ 
(i-r)(i-r+\gamma')\vert L_{ir} \vert \vert R_i \vert}  \label{eqn-bounds2}
\end{equation}
The following lemmas are needed to find a bound to (\ref{eqn-bounds2}).
\begin{lem}\label{lem-conv2}
Let $N=\gamma'-1/2$ and $i>2N+1$. Then, define $j^*=\lfloor \frac{i-1+N}{2} \rfloor$. It holds that $\vert L_{ij}  \vert \leq \vert L_{ij^*} \vert$.
\end{lem}
\begin{proof}
Consider the ratio $\frac{\vert L_{ij+1}  \vert}{\vert L_{ij}  \vert}$. It is easy to see that
\begin{equation} \frac{\vert L_{ij+1}  \vert}{\vert L_{ij}  \vert} = \frac{(i-j+\gamma')(i-j)}{\vert j+1-\gamma' \vert (j+1)} \label{eq-ratio}
\end{equation}
Now if $j \leq N < i/2-1/2 $, then $\vert j+1-\gamma' \vert=\gamma' -j -1>\gamma'-1=N+1/2$. One has then
\begin{equation} \frac{\vert L_{ij+1}  \vert}{\vert L_{ij}  \vert} = \frac{(i-j+\gamma')(i-j)}{(\gamma' -j -1) (j+1)} > \frac{i/2 (i/2 +\gamma')}{(N+1/2) (N+1)} \geq 1
\end{equation}
Thus the sequence always increases as long as $j\leq N$, and we can look for a maximum $j^*>N$. Then, for $j>N$, denote the ratio of (\ref{eq-ratio}) by $f$:
\begin{equation} f=\frac{\vert L_{ij+1}  \vert}{\vert L_{ij}  \vert} = \frac{(i-j+\gamma')(i-j)}{( j+1-\gamma' ) (j+1)} \label{eqn-ratioL}
\end{equation}
Now, $f\leq 1$ implies $(i-j)^2+\gamma' (i-j) \leq (j+1)^2 -\gamma' (j+1)$. Thus, $(i-j)^2-(j+1)^2 +\gamma' (i+1) \leq 0$. Manipulating the expression, one finds
$(i^2-1)-2j(i+1) +\gamma' (i+1) \leq 0$ and cancelling the term $(i+1)$ the following inequality for $j$ is reached:
\begin{equation}
j\leq \frac{i-1+\gamma'}{2} \label{eqn-bnd}
\end{equation}
Therefore, if (and only if) the bound given by (\ref{eqn-bnd}) is verified, $\frac{\vert L_{ij+1}  \vert}{\vert L_{ij}  \vert} \leq 1$. Therefore one concludes that the maximum of the sequence $\vert L_{ij} \vert$ is reached at
\begin{equation}
j=j^*=\lfloor \frac{i-1+\gamma'}{2}\rfloor =\lfloor \frac{i-1+N}{2}+1/4\rfloor=\lfloor \frac{i-1+N}{2}\rfloor  \label{eqn-bnd2}
\end{equation}
thus finishing the proof.
\end{proof}
\begin{lem}\label{lem-conv3}
Let $N=\gamma'-1/2$ and $i>2N+1$. Then,  one has that  $\frac{   \sum_{j=r+1}^{i} \vert L_{ij} \vert   \sum_{j=0}^r  \vert L_{ij} \vert}{(i-r)(i-r+\gamma')\vert L_{ir}} \leq 2 L_{ij^*}$, where $j^*$ is defined in (\ref{eqn-bnd2}).
\end{lem}
\begin{proof}

Now, to bound the term  $\sum_{j=r+1}^{i} \vert L_{ij} \vert   \sum_{j=0}^r  \vert L_{ij}$, consider two possibilities, and use Lemma~\ref{lem-conv2}. If $r < j^*$, then  $\sum_{j=r+1}^{i} \vert L_{ij} \vert   \sum_{j=0}^r  \vert L_{ij} \leq (i-r)(r+1) \vert L_{ij^*} \vert \vert L_{ir} \vert $. On the other hand, if $r \geq j^*$, then $\sum_{j=r+1}^{i} \vert L_{ij} \vert   \sum_{j=0}^r  \vert L_{ij} \leq (i-r)(r+1) \vert L_{ir+1} \vert  \vert L_{ij^*} \vert $. Therefore, if $r<j^*$:
$$
   \frac{   \sum_{j=r+1}^{i} \vert L_{ij} \vert   \sum_{j=0}^r  \vert L_{ij} \vert}{(i-r)(i-r+\gamma') \vert L_{ir} \vert}  
   \leq    \frac{ (r+1) \vert L_{ij^*} \vert }{i-r+\gamma'}  
      \leq    \frac{ (j^*+1) \vert L_{ij^*} \vert }{i-j^*+\gamma'}  
$$
else, if $r \geq j^*$,
$$
   \frac{   \sum_{j=r+1}^{i} \vert L_{ij} \vert   \sum_{j=0}^r  \vert L_{ij} \vert}{(i-r)(i-r+\gamma')\vert L_{ir} \vert}  
   \leq  \frac{ (r+1)  \vert L_{ir+1} \vert   \vert L_{ij^*} \vert }{(i-r+\gamma')\vert L_{ij} \vert}  
$$
and using (\ref{eqn-ratioL}),
$$
   \frac{   \sum_{j=r+1}^{i} \vert L_{ij} \vert   \sum_{j=0}^r  \vert L_{ij} \vert}{(i-r)(i-r+\gamma')\vert L_{ir} \vert}  
   \leq  \frac{ (i-r)   \vert L_{ij^*} \vert }{( r+1-\gamma' )  }  \leq  \frac{ (i-j^*)   \vert L_{ij^*} \vert }{( j^*+1-\gamma' ) } 
$$
Now, since $\frac{i-1+N}{2} \leq j^* \leq \frac{i+N}{2}$, one has that 
$$ \frac{ j^*+1 }{i-j^*+\gamma'} \leq \frac{\frac{i+N}{2}+1}{ i-\frac{i+N}{2}+N+1/2}
=\frac{i+N+2}{i+N+1}<2,
$$
and similarly,
$$
 \frac{i-j^*}{ j^*+1-\gamma'} \leq  \frac{i-\frac{i-1+N}{2}}{\frac{i-1+N}{2}+1/2-N}
 =\frac{i+1-N}{i-N}<2,
 $$
 thus concluding the proof.
\end{proof}
Thus we are left with showing that  $ \frac{ \vert L_{ij^*} \vert}{\vert R_i \vert}$ is decreasing, which is expected, since $L_{ij^*}$ is one the elements appearing in the sum $R_i$. Using the expression (\ref{eqn-explicitRi}) and the formula for $L_{ij}$ involving the Gamma function, one reaches:
\begin{eqnarray}
 \frac{ \vert L_{ij^*} \vert}{\vert R_i \vert}
 &=& \frac{i!^2  \Gamma(1-\gamma')\Gamma(i+1+\gamma')\Gamma(i+1-\gamma')}{2i! j^*! (i-j^*)! \Gamma(\gamma'-1) \Gamma(j^*+1-\gamma')\Gamma(i-j^*+1+\gamma')}.
\end{eqnarray}
Now, the decreasing character of the sequence is established as follows (consider the case where $i+N$ is odd so that $j^*=\frac{i-1+N}{2}$; the even case is analogous). Consider Stirling's approximation to the factorial, namely $n!\approx \sqrt{2 \pi n} \left(\frac{n}{\mathrm{e}} \right)^n$.
Then:
\begin{equation}
\frac{i!^2}{2i! j^*! (i-j^*)!}=
\frac{i!^2}{2i! \left( \frac{i-1+N}{2} \right)! \left(\frac{i+1-N}{2}\right)!}
 \approx 
 \sqrt{\frac{i}{\pi (i^2-(N-1)^2) }} \frac{1}{2^{2i}} \frac{\mathrm{e}^i}{ \left(\frac{i-1+N}{2} \right)^{\frac{i-1+N}{2} } \left(\frac{i+1-N}{2}\right)^{\frac{i+1-N}{2}}}.
 \label{eqn-stirling1}
\end{equation}
On the other hand, Stirling's approximation to the  Gamma function~\cite[p.257]{abramowitz} reads $\Gamma(z)\approx \sqrt{\frac{2\pi}{z}}\left( \frac{z}{\mathrm{e}} \right)^z$,
thus
\begin{eqnarray}
 \frac{\Gamma(i+1-\gamma')\Gamma(i+1+\gamma')}{\Gamma(j^*+1-\gamma')\Gamma(i-j^*+1+\gamma')}
 &=&\frac{\Gamma(i+1/2-N)\Gamma(i+3/2+N)}{\Gamma(\frac{i-N}{2})\Gamma(\frac{i+N}{2}+2)} \nonumber \\
 &\approx&  \frac{1}{\mathrm{e}^i}
  \sqrt{\frac{(\frac{i-N}{2})(\frac{i+N}{2}+2)}{(i+1/2-N)(i+3/2+N)}} \nonumber \\ && \times
  \frac{(i+1/2-N)^{i+1/2-N} (i+3/2+N)^{i+3/2+N}}{(\frac{i-N}{2})^{\frac{i-N}{2}} (\frac{i+N}{2}+2)^{\frac{i+N}{2}+2}} .
  \label{eqn-stirling2}
\end{eqnarray}
Putting together (\ref{eqn-stirling1})--(\ref{eqn-stirling2}), we obtain
\begin{eqnarray}
 \frac{ \vert L_{ij^*} \vert}{\vert R_i \vert}
 &\approx& \frac{i\Gamma(1-\gamma')}{ \sqrt{\pi} \Gamma(\gamma'-1)} \sqrt{f_1(i)} f_2(i)^i f_3(i),
\end{eqnarray}
where we have broken the approximation in three functions:
\begin{eqnarray}
f_1(i)&=&\frac{i}{\pi (i^2-(N-1)^2) } \frac{(\frac{i-N}{2})(\frac{i+N}{2}+2)}{(i+1/2-N)(i+3/2+N)},\\
f_2(i)&=& \frac{(i+1/2-N)(i+3/2+N)}{\sqrt{(i-1+N)(i+1-N)( i-N) (i+N+4)}},\\
f_3(i)&=& \frac{(i+1/2-N)^{1/2-N} (i+3/2+N)^{3/2+N}}{ \left(\frac{i-1+N}{2} \right)^{\frac{-1+N}{2} } \left(\frac{i+1-N}{2}\right)^{\frac{1-N}{2}}
(\frac{i-N}{2})^{\frac{-N}{2}} (\frac{i+N}{2}+2)^{\frac{N}{2}+2}
}.
\end{eqnarray}
Notice that, clearly, $
\lim_{i \rightarrow \infty} f_1(i)=0$ (since $f_1(i)$ behaves like $\mathcal{O}(1/i)$ for large $i$), 
$\lim_{i \rightarrow \infty} f_2(i)=1$, and 
$\lim_{i \rightarrow \infty} f_3(i)= 16$,
thus it only remains to compute $\lim_{i \rightarrow \infty} f_2(i) ^i $, which is an indeterminate of the kind $1^\infty$. Resolving it (the details are omitted for brevity) one obtains that the limit is indeed 1. Thus, it is possible to find the decreasing sequence $c_i$ in (\ref{eq-bounds}), concluding the proof of convergence and analyticity in odd dimension.

\subsubsection{Proof of analyticity for even dimension}
The fact that $\gamma'$ is integer makes the odd approach a priori impossible since (\ref{def-Lij}) would not be well-defined (it would contain divisions by zero). However, to overcome that difficulty we employ a partial solution for the kernel equations, to order $\gamma'-1$, that helps regularize the problem.

For that end, consider $F(r,\rho)=\sum_{i=0}^{\gamma'-1} r^{2i} \phi_i(\rho^2)$. Replacing this function in (\ref{eqn-gn})--(\ref{eqn-Gnbc}) results in
$$
\sum_{i=0}^{\gamma'-1} r^{2i} \frac{\lambda(\rho^2)+c}{\epsilon} \phi_i(\rho^2) =  \sum_{i=0}^{\gamma'-1}
\left[(4i(i-\gamma')   r^{2i-2} \phi_i(\rho^2) -
r^{2i} \left(  4 \rho^2 \phi_i'' (\rho^2) +2 (2+\gamma') \phi' (\rho^2) \right) \right]
$$
and one gets the following recursive set of ODEs:
$$
\frac{\lambda(\rho^2)+c}{\epsilon} \phi_i(\rho^2) 
= 4(i+1)(i+1-\gamma')  \phi_{i+1}(\rho^2) -4 \rho^2 \phi_i'' (\rho^2) - 2 (2+\gamma') \phi'_{\gamma'-1} (\rho^2) $$
which is solved starting at $i=\gamma'-1$:
$$
\frac{\lambda(\rho^2)+c}{\epsilon} \phi_{\gamma'-1}(\rho^2) 
= -4 \rho^2 \phi_{\gamma'-1}'' (\rho^2) - 2 (2+\gamma') \phi'_{\gamma'-1} (\rho^2) $$
This can be written as
$$
4 x \phi_{\gamma'-1}'' +2(2+\gamma') \phi'_{\gamma'-1}+\frac{\lambda(x)+c}{\epsilon} \phi_{\gamma'-1}=0
$$
which is an ODE with a  regular singular point at $x=0$. By applying the Frobenius method~\cite[Chapter 36]{howell} one can rewrite this equation as
$$
4 x^2 \phi_{\gamma'-1}'' +2x(2+\gamma') \phi'_{\gamma'-1}+\frac{\lambda(x)+c}{\epsilon}x \phi_{\gamma'-1}=0
$$
and its indicial equation is $r(r-1)+(1+\gamma'/2)r=0$, thus $r_1=0$ and $r_2=\gamma'/2$ (non-integer). We are interested in the solution of the form $\phi_{\gamma'/2}=\sum_{i=0}^\infty a_i \rho^{2i}$ and discard the other solution. By Fuchs' theorem~\cite[p.146]{butkov} this solution is analytic where $\lambda(x)$ is analytic, thus the radius of convergence of the resulting $\phi_{\gamma'-1}(\rho^2)$ is greater than one.
Next, for $i=\gamma'-2$ up to $i=0$:
$$
4 x \phi_{i}'' +2(2+\gamma') \phi'_{i}+\frac{\lambda(x)+c}{\epsilon} \phi_{i}=4(i+1)(i+1-\gamma')  \phi_{i+1}(x)
$$
which, has the same indicial equation and again, also admits a solution in the required form. Applying once more Fuchs' theorem, this solution is analytic in intervals where both $\lambda(x)$ and $\phi_{i+1}$ are analytic. Thus by induction we find a family of solutions such that the radius of convergence of all the $\phi_i$ is greater than one.

The solutions just found have a degree of freedom (the first coefficient $a_i$ of their power series, which is $\phi_i(0)$). The idea is to construct the solution such that the boundary condition $G_{lm}^n(r,r)=H(r)$ is satisfied up to order $2\gamma'-2$. Thus:
$F(r,r)=\sum_{i=0}^{\gamma'-1} r^{2i} \phi_i(r^2)$ and expanding in power series $\phi_i(r^2)$:
$$
F(r,r)=\sum_{i=0}^{\gamma'-1} r^{2i}\left(
 \phi_i(0)+\frac{r^2}{1!} \phi_i'(0)+\frac{r^4}{2!} \phi_i''(0)+\hdots\right)$$
 Thus:
\begin{eqnarray}
\phi_0(0)&=&H(0),\\
\phi_1(0)+\frac{1}{1!} \phi'_0(0)&=&\frac{1}{1!} H'(0),\\
\phi_2(0)+\frac{1}{1!} \phi'_1(0)+\frac{1}{2!}\phi''_0(0)&=&\frac{1}{2!} H''(0),\\
&\hdots& \\
\phi_{\gamma-1} (0)+\hdots + \frac{1}{(\gamma-1)!} \phi^{(\gamma-1)}_0(0)&=&\frac{1}{(\gamma-1)!} H^{(\gamma-1)}(0).
\end{eqnarray}
It can be shown this scheme produces valid initial values for the $\phi_i$'s. However, an easier approach is following the general series approach of Section~\ref{sec-series} up to order $i=\gamma'-1$. By uniqueness of the series development and identifying coefficients, it can be easily shown that $\phi_i(0)=C_{ii}$.

Next, calling $G_{lm}^n (r,\rho)=\check G_{lm}^n (r,\rho)+F(r,\rho)$ the new boundary condition for the PDE becomes:
$\check G_{lm}^n(r,r)=H(r)-F(r,r)$ which starts at order $2\gamma'$. Thus one can propose $\check G_{lm}^n(r,\rho)=r^{\gamma} F_2(r,\rho)$. One can see the PDE for $F_2$ is
\begin{eqnarray}\label{eqn-gn2}
\frac{\lambda(\rho)+c}{\epsilon} F_2(r,\rho) &=& 
 \partial_{rr}F_2(r,\rho) +(1+\gamma) \frac{\partial_rF_2(r,\rho)}{r}
-\partial_{\rho\rho} F_2(r,\rho) -(1+\gamma) \frac{ \partial_\rho  F_2(r,\rho)}{\rho},
\end{eqnarray}
and following previous Sections, and calling $\psi(r^2)=\frac{H(r)-F(r,r)}{r^{2\gamma'}}$ and abusing the notation by keeping the same name for the coefficients $C_{ij}$, one can find a power series development $F_2(r,\rho)=\sum_{i=0}^\infty \left(\sum_{j=0}^i C_{ij} r^{2j} \rho^{2(i-j)} \right)$ as
\begin{eqnarray}
\forall i,\quad && \sum_{j=0}^i C_{ij} =- \frac{\psi_i}{2(2i+1)}, \\
\forall i\geq1,0\leq j \leq i-1,&& \,\,  (j+1)(j+1+\gamma')C_{i(j+1)} -(i-j)(i-j+\gamma')C_{ij} =
B_{{(i-1})j},
\end{eqnarray}
Now the approach of Section~\ref{sec-odd} becomes applicable and even easier, since all coefficients are positive. Indeed, define
\begin{equation}
L_{i0}=1, \qquad 
L_{ij}=\left(\begin{array}{c} i \\ j\end{array}\right) \frac{(i+\gamma')(i-1+\gamma')\hdots(i-j+\gamma'+1)}{(1+\gamma')(2+\gamma')\hdots(j+\gamma')}>0, \quad j>0
\label{def-Lij2}
\end{equation}
Mimicking  Section~\ref{sec-odd} we reach
\begin{eqnarray}
\alpha_n&=& \sum_{j=0}^i \vert C_{ij}\vert \nonumber \\ &\leq &
 \frac{\vert \psi_i \vert }{2 (2i+1)} \frac{ \sum_{j=0}^i \vert L_{ij} \vert}{\vert R_i \vert} +\sum_{r=0}^{i-1}\vert B_{{(i-1})r} \vert \left(  \frac{  \vert R_{ir} \vert \sum_{j=r+1}^{i} \vert L_{ij} \vert + \vert R_i-R_{ir}\vert \sum_{j=0}^r  \vert L_{ij} \vert}{ 
(i-r)(i-r+\gamma')\vert L_{ir} \vert \vert R_i \vert }\right)
\nonumber \\ &= &
 \frac{\vert \psi_i \vert }{2 (2i+1)}+\sum_{r=0}^{i-1}\vert B_{{(i-1})r} \vert  \frac{  2  R_{ir} ( R_i-R_{ir})}{ 
(i-r)(i-r+\gamma') L_{ir}  R_i } \label{eqn-alphandef-even}
,
\end{eqnarray}
where the last step can be carried out due to positivity of the redefined coefficients $L_{ij}$ compared to Section~\ref{sec-odd}.
Again, we apply Lemma~\ref{lem-conv} to (\ref{eqn-alphandef-even}). In this case, we define 
$$
b_i=\frac{ 1}{2 \epsilon(2i+1)} \geq 0
$$
which is already a decreasing sequence. Then we need to find $c_i$ such that
\begin{equation}
 \max_{r\in\{0,\hdots,i-1\}} \frac{  2  R_{ir} ( R_i-R_{ir})}{ 
(i-r)(i-r+\gamma') L_{ir}  R_i}  \leq c_i \label{def-ci-even}
\end{equation}
and $c_i$ needs to be proven decreasing (for sufficiently large $i$) and convergent to zero. Consider the following lemma.
\begin{lem} \label{Lem_Lijodd}
Consider $L_{ij}$ as defined in (\ref{def-Lij2}) and $R_i$, $R_{ij}$ as defined in (\ref{def-Lij2}). Then:
\begin{enumerate}
\item $L_{i(i-j)}=L_{ij}$
\item $R_i-R_{ir}=R_{i(i-r-1)}$
\item Let $F(i,r)=\frac{  R_{ir} ( R_i-R_{ir})}{ 
(i-r)(i-r+\gamma')L_{ir}}$ for $i,r$ nonnegative integers with $r< i$. Then one has $F(i,r)= F(i,i-r-1)$.
\item If $r\leq i/2$, $R(i,r) \leq (r+1) L_{ir}$.
\item  It holds that $\max_{r\in\{0,\hdots,i-1\}} \frac{  2  R_{ir} ( R_i-R_{ir})}{ 
(i-r)(i-r+\gamma') L_{ir} } \leq   2\frac{i+2}{i(i+2\gamma')} R_i$.
\end{enumerate}
\end{lem}
\begin{proof}
Writing  $L_{ij}=\left(\begin{array}{c} i \\ j\end{array}\right)\frac{\Gamma(1+\gamma')\Gamma(i+1+\gamma')}{\Gamma(j+1+\gamma')\Gamma(i-j+1+\gamma')}$ the first property is evident, whereas the second property is immediate from the first one since $R_i-R_{ij}= \sum_{k=j+1}^{k=i}L_{ik}=\sum_{k=0}^{i-j-1} L_{i(i-k)}=\sum_{k=0}^{i-j-1} L_{ik}=R_{i(i-j-1)}$.

For the third property, note that 
\begin{eqnarray}
F(i,r)&=&\frac{(r+1)(r+1+\gamma)}{(i-r)(i-r+\gamma')}\frac{L_{i(i-r-1)}}{L_{ir}}F(i,i-r-1) \nonumber \\
&=&\frac{(r+1)(r+1+\gamma)}{(i-r)(i-r+\gamma')}\frac{L_{i(r+1)}}{L_{ir}}F(i,i-r-1)
 \nonumber \\
&=&F(i,i-r-1)
\end{eqnarray}
The fourth property is obvious noting that $L_{ij} \leq L_{i(j+1)}$ for $j <  i/2 $. Finally, for the last property, first note that it is only required to study $0 \leq r \leq   i/2  $ given the third property. Now:
\begin{eqnarray}
F(i,r)&=&\frac{  R_{ir} ( R_i-R_{ir})}{ 
(i-r)(i-r+\gamma')L_{ir}} \nonumber \\
&\leq &
\frac{  r+1 R_{i(i-r-1)}}{ 
(i-r)(i-r+\gamma')}
\nonumber \\
&\leq &
\frac{r+1}{ (i-r)(i-r+\gamma'} R_i
\end{eqnarray}
and since this is an increasing function of $r$ for $0\leq r<i$, we can bound it by its value at $r=i/2$, thus proving the final property.
\end{proof}
Therefore, setting $c_i=4\frac{i+2}{i(i+2\gamma')}$, a sequence decreasing to zero, we can apply Lemma~\ref{lem-conv} to (\ref{eqn-alphandef-even}) and follow the same steps as in Section~\ref{sec-odd} to obtain the result of Theorem~\ref{th-Gresult} for the even dimension case.
\section{Observer design}\label{sec-observer}
This section designs an observer for \eqref{eqn-un}-\eqref{eqn-un-bc} from the measured output $\partial_r u^{m}_{l}(t,R)$ as follows:
\begin{align}\label{eqn-ul-obmain}
\partial_t \hat{u}^{m}_{l}=&\epsilon \frac{\partial_r \left( r^{n-1} \partial_r \hat{u}^{m}_{l} \right)}{r^{n-1}} - l(l+n-2) \frac{\epsilon}{r^2} \hat{u}^{m}_{l}+\lambda(r)\hat{u}^{m}_{l} + p_{lm}^n(r)(\partial_r u^{m}_{l}(t,R)-\partial_r \hat{u}^{m}_{l}(t,R)),
\end{align}
with boundary condition
\begin{equation}\label{eqn-ul-obbc}
\hat{U}^{m}_{l}(t)=U^{m}_{l}(t).
\end{equation}
We need to design the output injection gain $p_{lm}^n(r)$. Closely following~\cite{nball}, define the observer error as $\tilde u=u-\hat u$. The observer error dynamics are given by
\begin{eqnarray}\label{eqn-utilden}
\frac{\partial \tilde u_{lmt}}{\partial t} &=&\frac{\epsilon}{r^{n-1}} \partial_r \left( r^{n-1}\partial_r\tilde  u_{lm} \right)-l(l+n-2) \frac{\epsilon}{r^2} \tilde  u_{lm}+\lambda(r) \tilde  u_{lm}-p_{lm}^n(r) \partial_r\tilde u_{lm}(t,R),
\end{eqnarray}
with  boundary conditions 
\begin{eqnarray}
\tilde u_{lm}(t,R)&=&0.
\end{eqnarray}
 Next we use the backstepping method to find a value of $p_{lm}^n(r)$ that guarantees convergence of $\tilde u$ to zero. This ensures that the observer estimates tend to the true state values. Our approach to design $p(r)$ is to seek a mapping that transforms (\ref{eqn-utilden}) into the following target system
\begin{eqnarray}\label{eqn-wtilden}
\frac{\partial \tilde w_{lmt}}{\partial t} &=&\frac{\epsilon}{r^{n-1}} \partial_r \left( r^{n-1} \partial_r \tilde w_{lm} \right)
-c \tilde w_{lm} -l(l+n-2) \frac{\epsilon}{r^2} \tilde w_{lm},
\end{eqnarray}
with boundary conditions
\begin{eqnarray}
\tilde w_{lm}(t,R)&=&0.
\end{eqnarray}

The transformation is defined as follows:
\begin{equation}\label{eqn-obstrans}
\tilde u_{lm}(t,r)=\tilde w_{lm}(t,r)-\int_r^RP^n_{lm}(r,\rho) \tilde w_{lm}(t,\rho) d\rho,
\end{equation}
and then $p_{lm}^n(r)$ will be found from transformation kernel as an additional condition.

From~\cite{nball}, one obtains the following PDE that the kernel must verify:
\begin{equation}\label{obs-kernel}
\frac{1}{r^{n-1}} \partial_r \left(r^{n-1}   \partial_r P^n_{lm} \right)
-  \partial_\rho \left(\rho^{n-1} \partial_\rho \left( \frac{P^n_{lm}}{\rho^{n-1}} \right)\right)
-l(l+n-2) \left(\frac{1}{r^2}- \frac{ 1}{\rho^2} \right)P^n_{lm}=-\frac{\lambda(r)}{\epsilon} P^n_{lm}
\end{equation}
In addition we find a value for the output injection gain kernel
\begin{equation}
p_{lm}^n(r)=\epsilon P^n_{lm}(r,R)  
\end{equation}

Also, the following boundary condition has to be verified
\begin{eqnarray}
0&=& \lambda(r)  +\epsilon \left( \partial_r   P^n_{lm} (r,\rho)\right) \bigg|_{\rho=r} +\frac{ \epsilon}{r^{n-1}}  \frac{d}{dr} \left(r^{n-1}   P^n_{lm}(r,r) \right) 
+\epsilon \partial_\rho
\left( \frac{P^n_{lm} (r,\rho)}{\rho^{n-1}}\right) \bigg|_{\rho=r} r^{n-1} ,
\end{eqnarray}
which can be written as
\begin{equation}
0= \lambda(r)  +  \epsilon  \partial_r P^n_{lm} (r,r) +\epsilon \frac{d}{dr}\left(P^n_{lm}(r,r)\right) +(n-1)\frac{ \epsilon P^n_{lm}(r,r)  }{r}  +\epsilon \partial_\rho P^n_{lm}(r,r) - (n-1)\frac{ \epsilon P^n_{lm}(r,r)  }{r}.\quad
\end{equation}
Following~\cite{nball}, and after some computations, we reach boundary conditions for the kernel equations as follows
\begin{eqnarray}\label{PBC-1}
P^n_{lm}(0,\rho)&=&0, \quad \forall l \neq 0\\\label{PBC-2}
P^n_{lmr} (0,\rho)&=&0 , \quad \forall l \neq 1\\\label{PBC-3}
P^n_{lm}(r,r) 
&=& -\frac{\int_0^r \lambda(\sigma) d\sigma  }{2 \epsilon }.
\end{eqnarray}
It turns out the observer kernel equation can be transformed into the control kernel equation, therefore obtaining a similar explicit result. For this, define
\begin{equation}
\check{P}^n_{lm}(r,\rho)=\frac{\rho^{n-1}}{r^{n-1}} P^n_{lm}(\rho,r),
\end{equation}
and it can be verified that the equation now governing $\check{P}^n_{lm}(r,\rho)$ is exactly  the equation satisfied by $K^n_{lm}(r,\rho)$. Thus $\check{P}^n_{lm}(r,\rho)=K^n_{lm}(r,\rho)$ and we can apply our previous result of Section~\ref{sec-wp}.

The observer error $\tilde{u}$ has the same stability properties derived in Section~\ref{sec-stab} for the closed-loop system under the full-state control. As in the controller case, only a limited number of modes need to be estimated; namely, those that are not naturally stable by virtue of Lemma~\ref{Lem_wl}, this being the main difference with the result given in~\cite{nball}.

Finally, the controller-observer augmented system can be proved closed-loop stable as in~\cite{nball}, using the separation principle given the linearity of the system, with desired convergence rate, and without much modification; we skip the details, which requires going up to $H^1$ stability, as in~\cite{nball}.

\section{Simulation Study}\label{simulation}
In this section, the simulation experiment on three-dimensional ball (n=3) is taken as an example to illustrate the effectiveness of proposed control. 
The system with the output-feedback  control law is simulated over $ 0\leq t\leq 2s$ with the following parameters: $\epsilon=1, \lambda(r)=10r^4+50r^2+50$, $c=3$. We consider that the system is initially at the random quantity, $u_0 \in [0,10]$, and the observer's initial condition is set as actual state plus an error of normal distribution with zero mean and $\sigma^2=0.5$.

\begin{figure} 
	\centering
	\footnotesize
	\begin{tabular}{@{}c@{}c@{}c@{}}
		\includegraphics[width=0.33\textwidth]{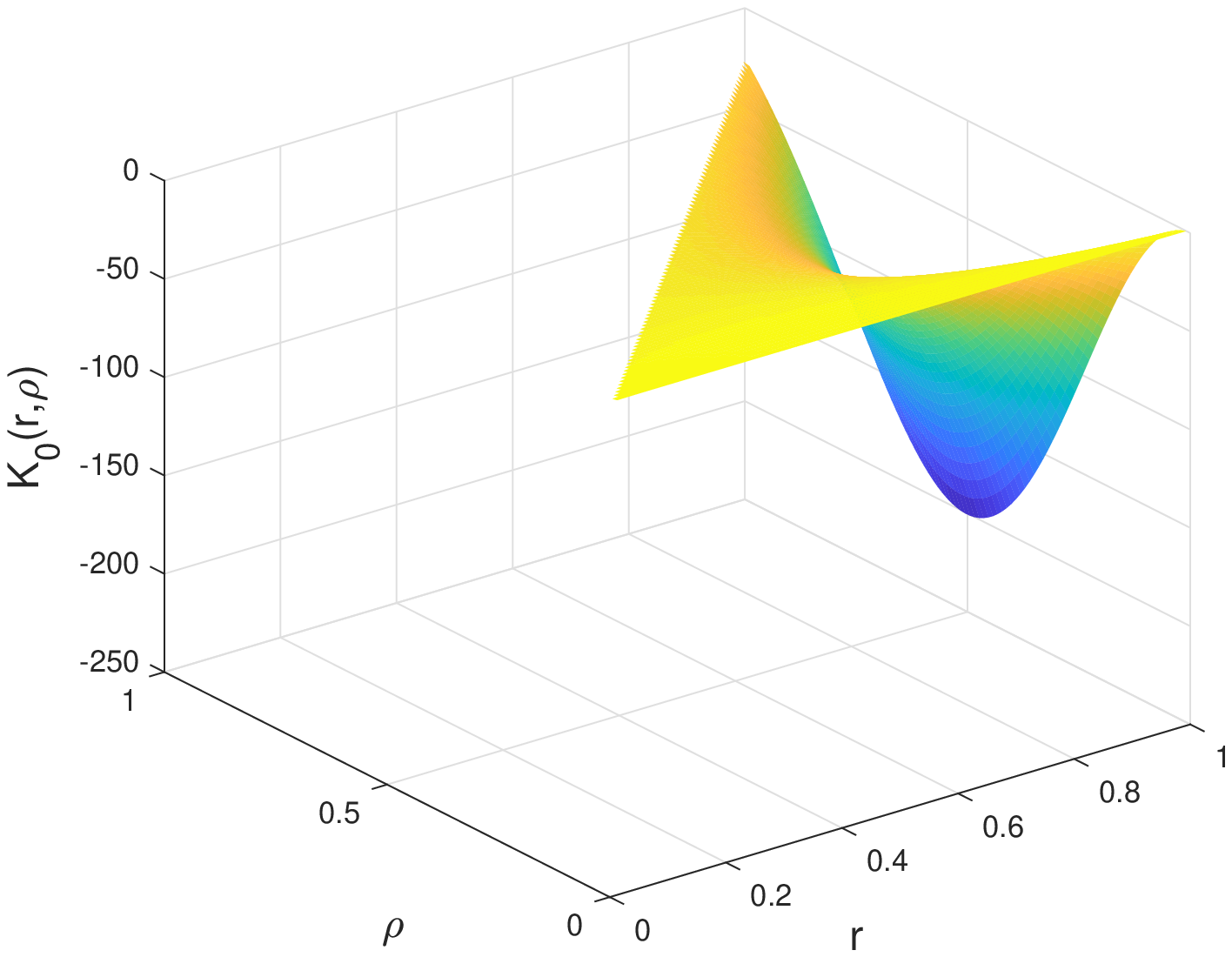}
		&\includegraphics[width=0.33\textwidth]{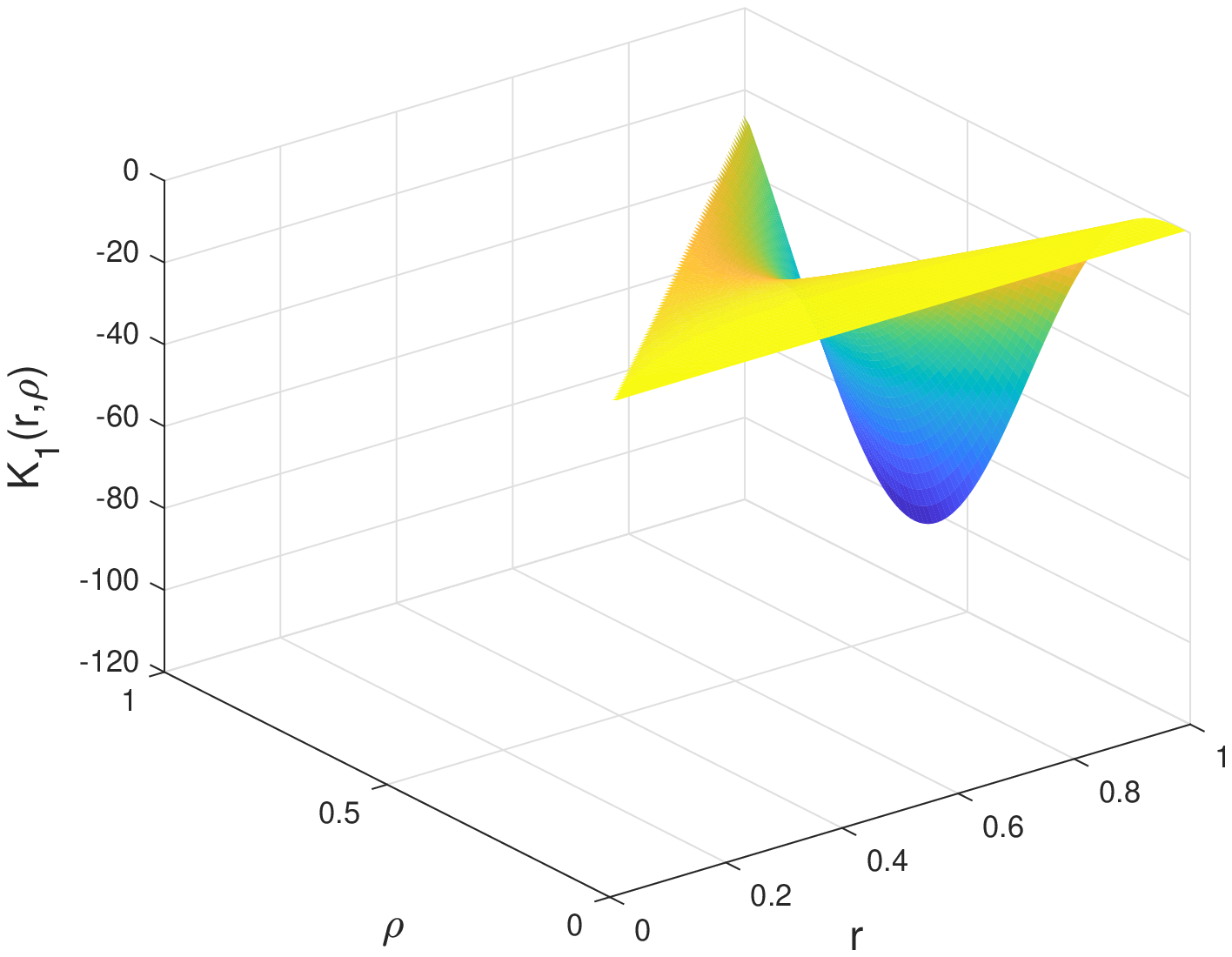}
		&\includegraphics[width=0.33\textwidth]{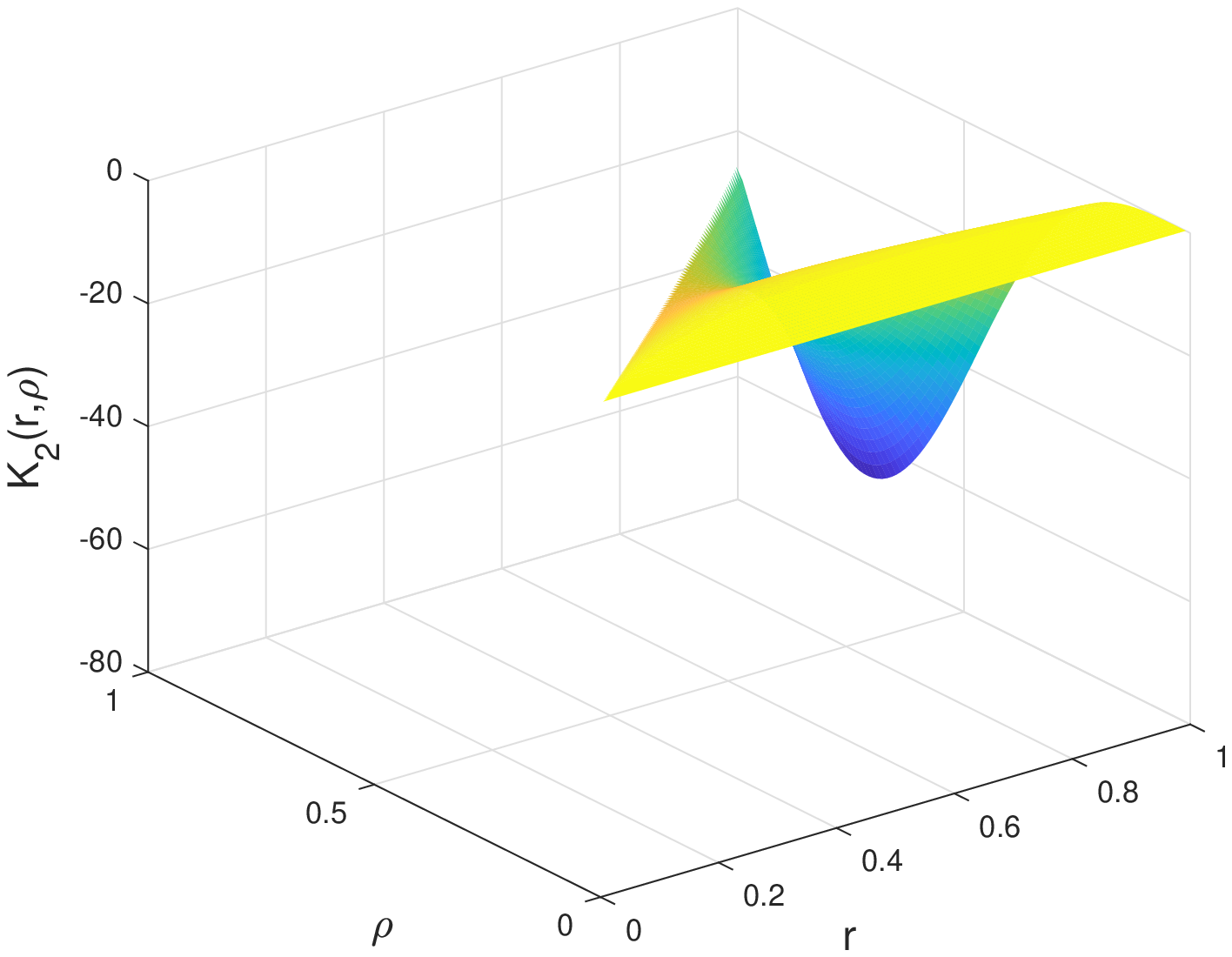}			\\
		(a) l=0 & (b) l=1 & (c) l=2 
	\end{tabular}\\
	\begin{tabular}{@{}c@{}c@{}c@{}}
		\includegraphics[width=0.33\textwidth]{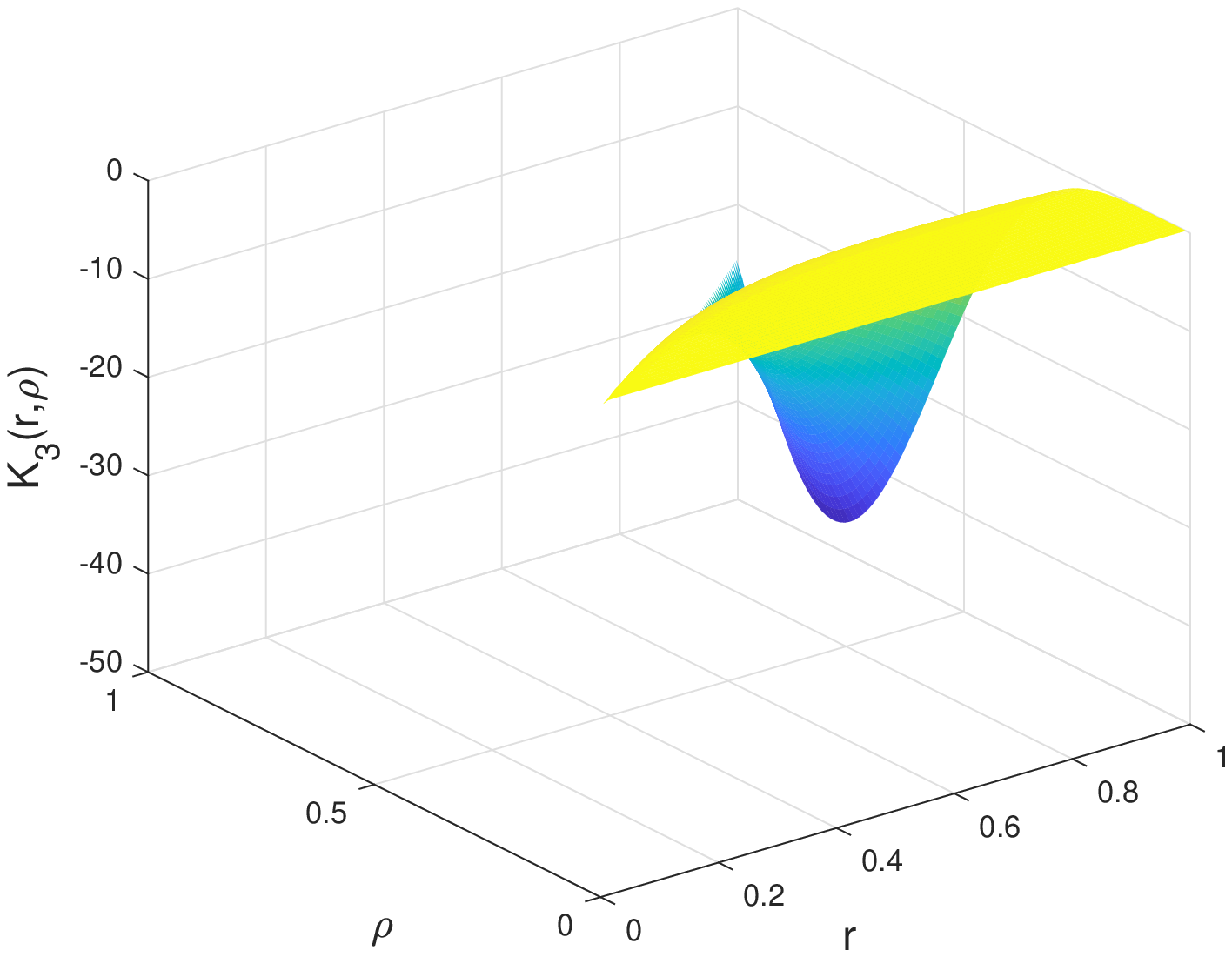}
		&\includegraphics[width=0.33\textwidth]{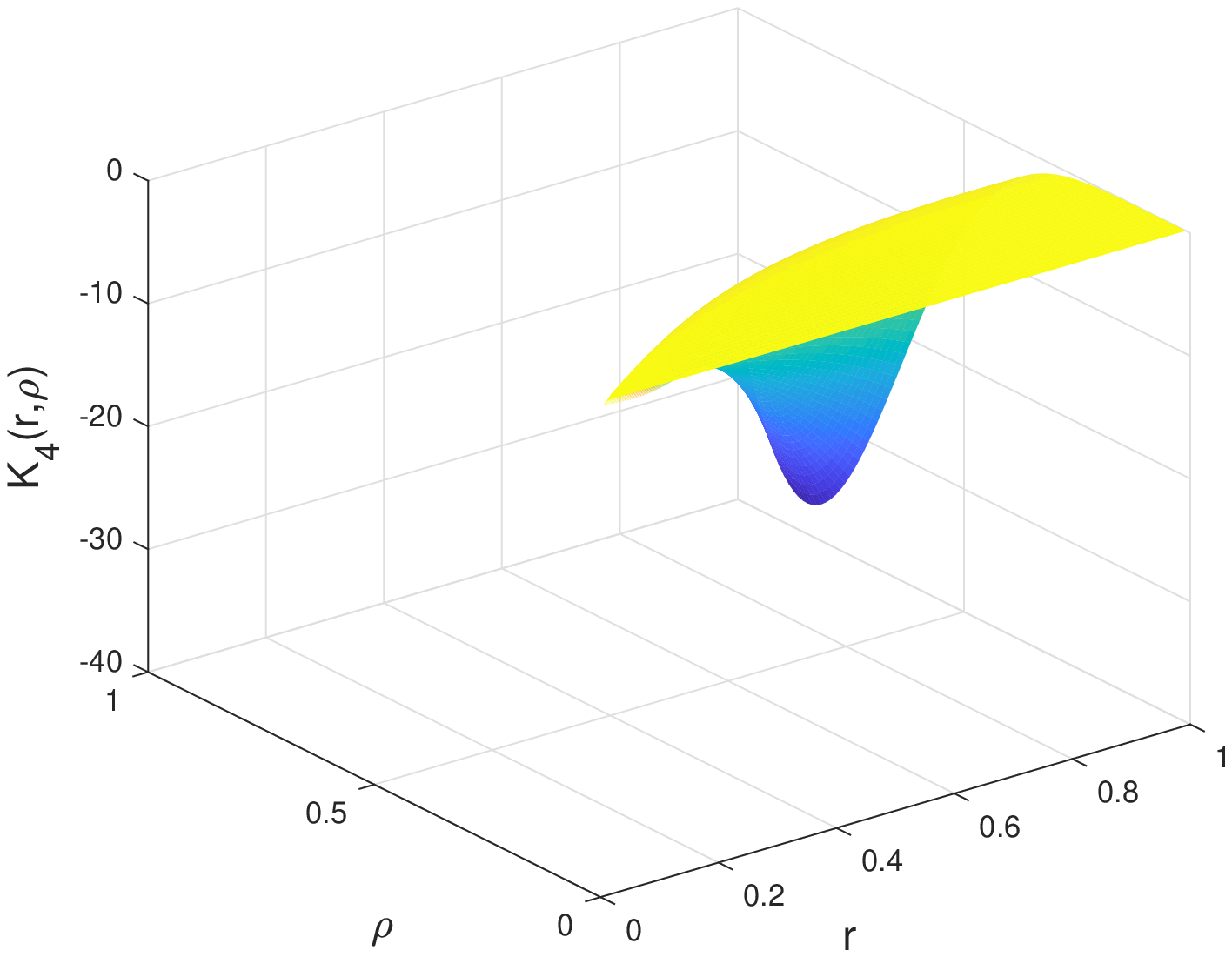}
		&\includegraphics[width=0.33\textwidth]{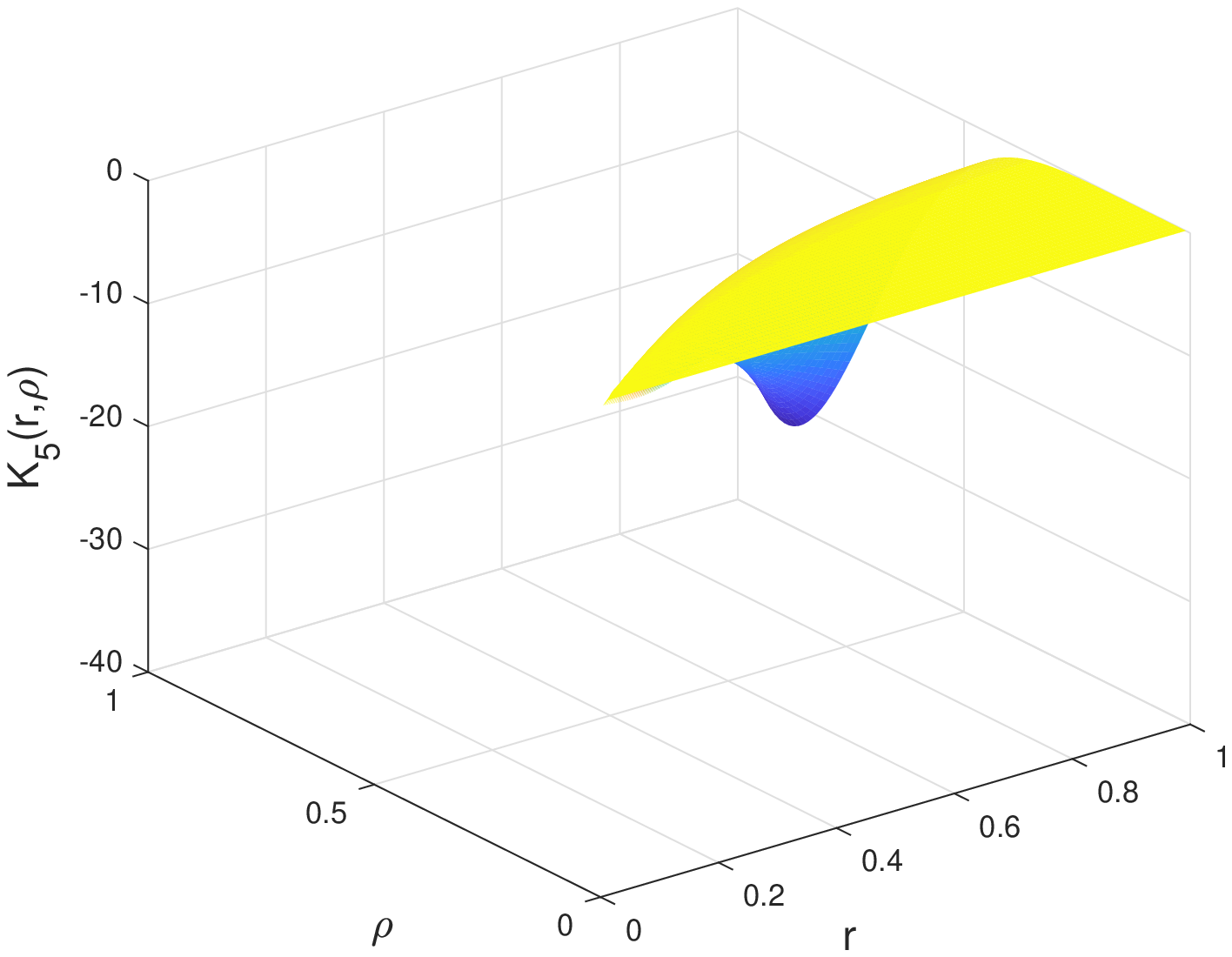}		\\ 
		(d) l=3 & (e) l=4 & (f) l=5
	\end{tabular}\\
	\caption{Polynomial approximation of the first six  control gains $K_l(r,\rho)$, $l=0,1,\dots,5$.}
	\label{fig_kernel}
\end{figure}

Fig.~\ref{fig_kernel} shows the plots of the polynomial approximation of kernels $K_{lm}^{3}$, which is obtained by using (\ref{def-power}) up to a cut-off at the $p$-th powers. The value of $K$ does not depend on $m$ so we omit that sub-index. The value of $p$ is chosen as $p=15$. Applying Lemma~\ref{lem_ul_open}, one can obtain $l$ to be $11$; however, here to save space, we only show the first six approximate numerical solutions of control gains. As shown in Fig.~\ref{fig_kernel}, we find that  the $K_l$ becomes increasingly smaller when $l$ increases.

In order to avoid a dramatic increase in the complexity of simulation caused by the high dimension, in our simulations we employ a method also based on spherical harmonics expansions which greatly reduces the error.
Thus, we only calculate the harmonics $u_{lm}$ which only need discretization in the radial direction, and then sum a finite number $S$ of harmonics to recover $u$.
When $S>0$ is a large enough integer, the error caused by the use of a finite number of harmonics is much smaller than the angular discretization error. Thus, the simulation is carried out using the formula
\begin{align}
u(t,r,\theta_1,\theta_2)=&\sum_{l=0}^{l=S} \sum_{m=-l}^{m=l} u_l^m(t,r) Y_{lm}^{3} (\theta_1,\theta_2)
\end{align}
where the spherical harmonics are defined as
\begin{equation}
Y_{lm}^{3}(\theta_1,\theta_2)=(-1)^m \sqrt{\frac{2l+1}{4\pi}\frac{(l-m)!}{(l+m)!}} P_{lm}^3(\cos(\theta_1)) e^{jm\theta_2}
\end{equation}
with $ P_{lm}^3$ the associated Legendre polynomial defined
as
\begin{equation}
P_{lm}^3(s)=\frac{1}{2^l l!}(1-s^2)^{m/2}\frac{d^{l+m}}{ds^{l+m}} (s^2-1)^l
\end{equation}

\begin{figure} 
	\centering
	\footnotesize
	\begin{tabular}{@{}c@{}c@{}}
		\includegraphics[width=0.5\textwidth]{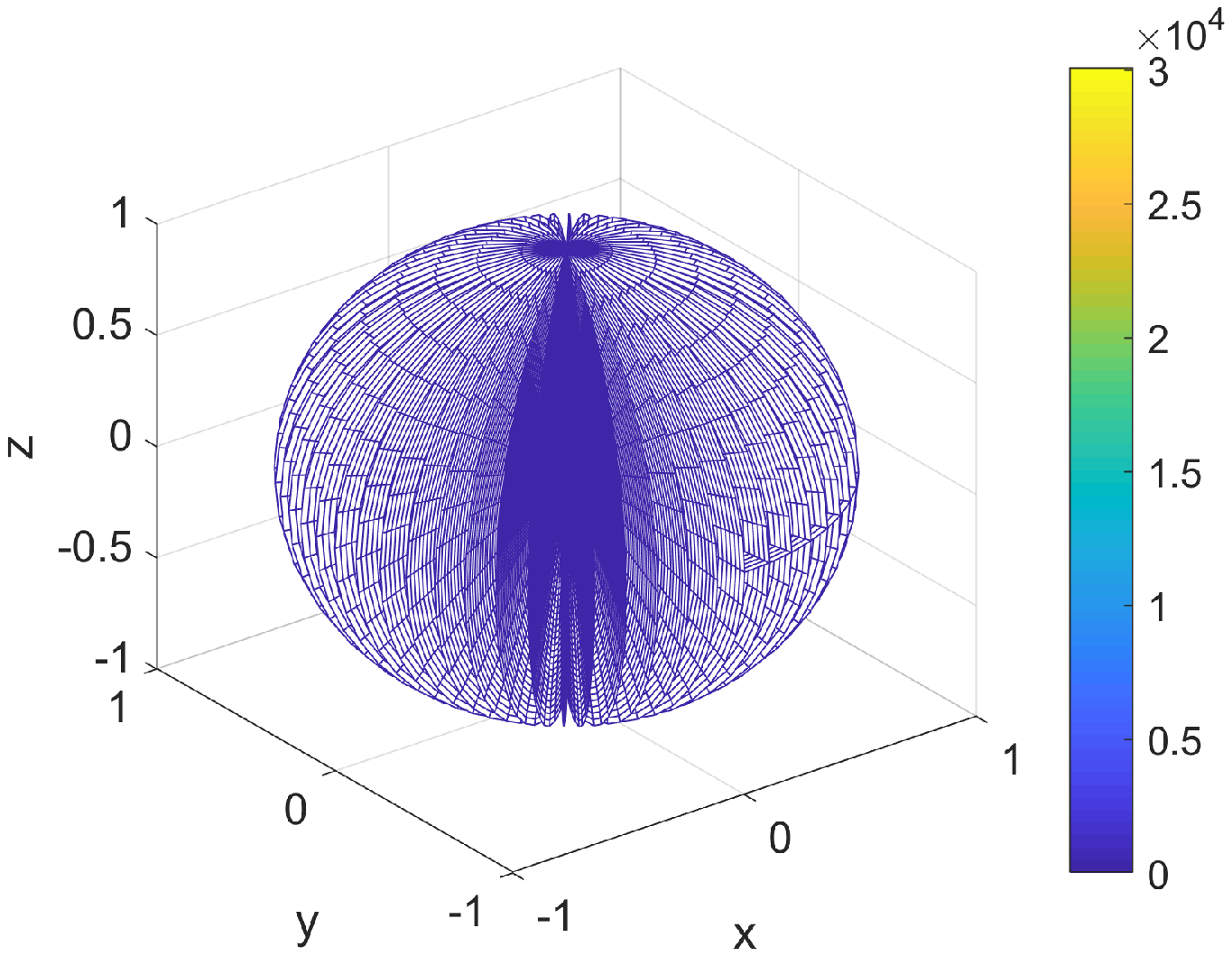}
		&\includegraphics[width=0.5\textwidth]{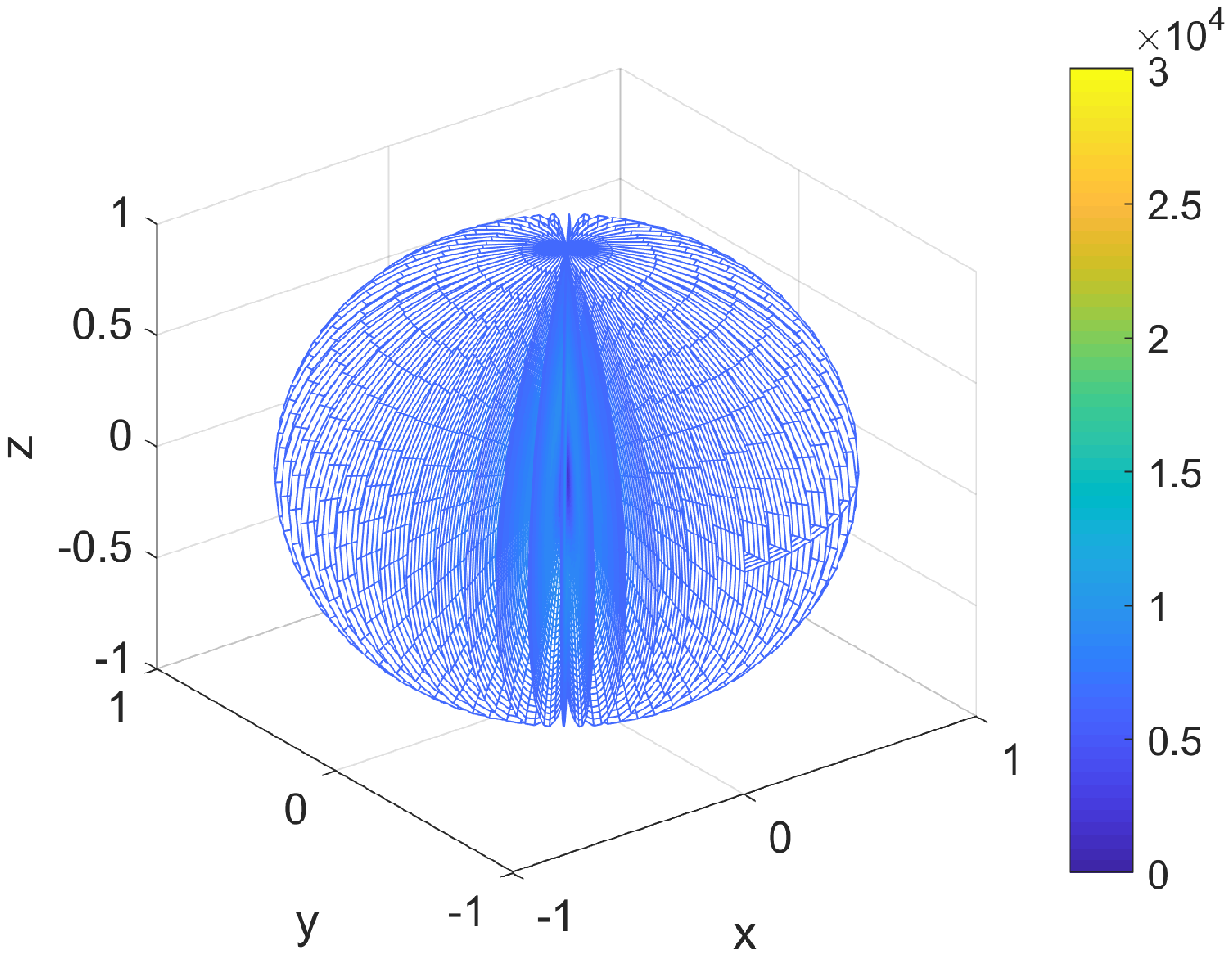}\\
		(a) t=0s & (b) t=0.18s
	\end{tabular}\\
	\begin{tabular}{@{}c@{}c@{}}
		\includegraphics[width=0.5\textwidth]{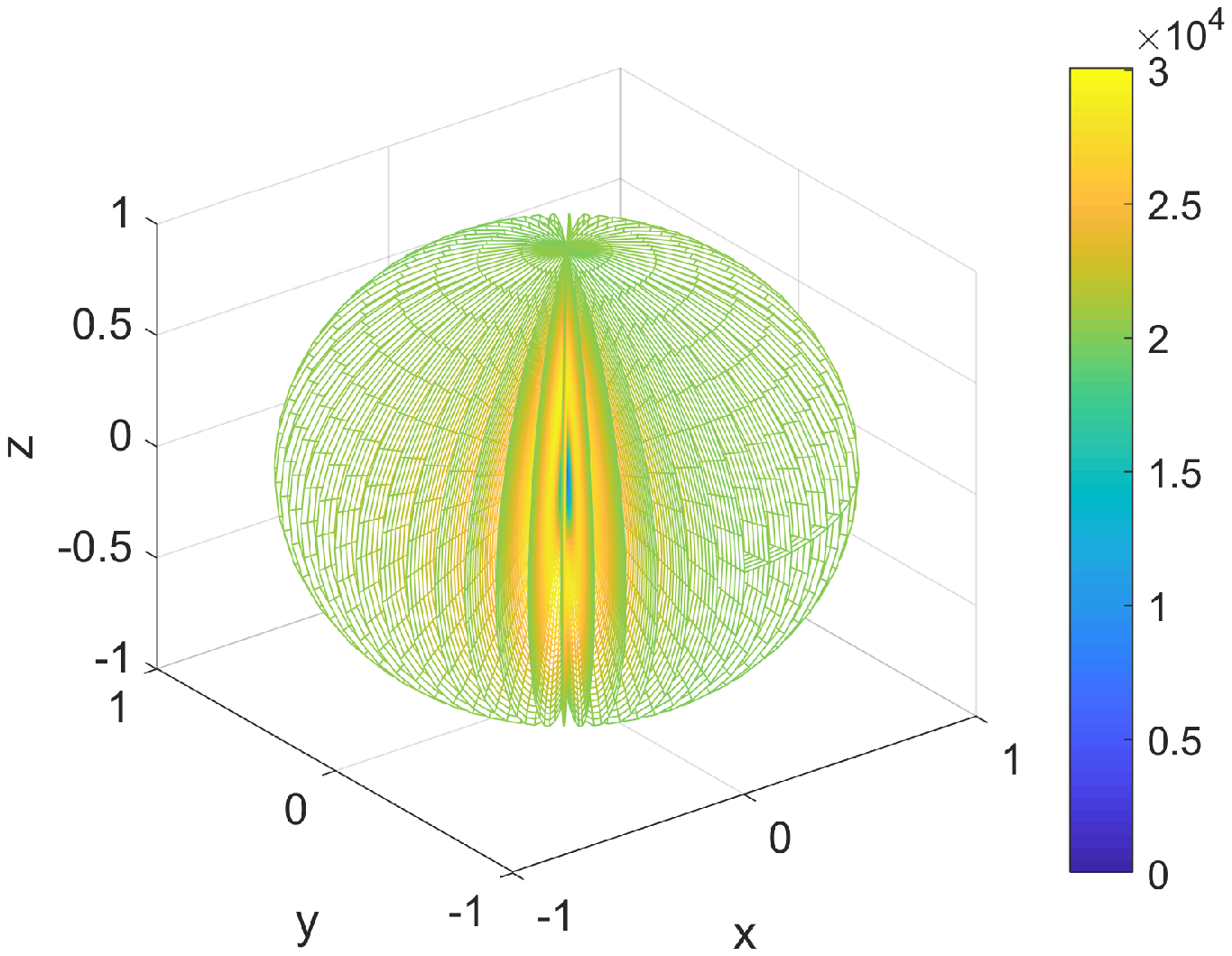}
		&\includegraphics[width=0.5\textwidth]{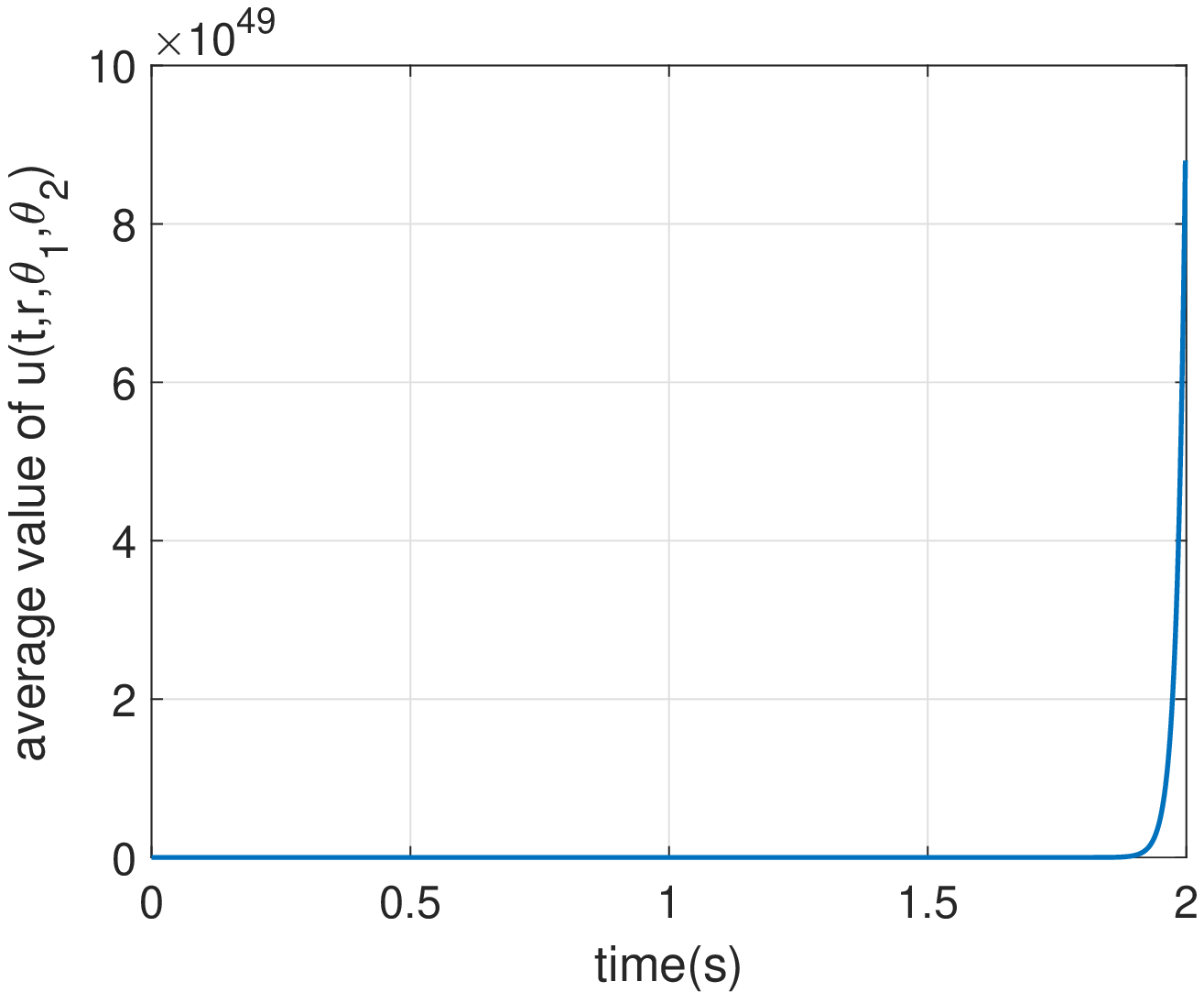}\\
		(c) t=0.2s & (d) Mean of the open-loop system
	\end{tabular}
	\caption{Open-loop evolution.}
	\label{fig_openloop}
\end{figure}

\begin{figure} 
	\centering
	\footnotesize
	\begin{tabular}{@{}c@{}c@{}}
		\includegraphics[width=0.5\textwidth]{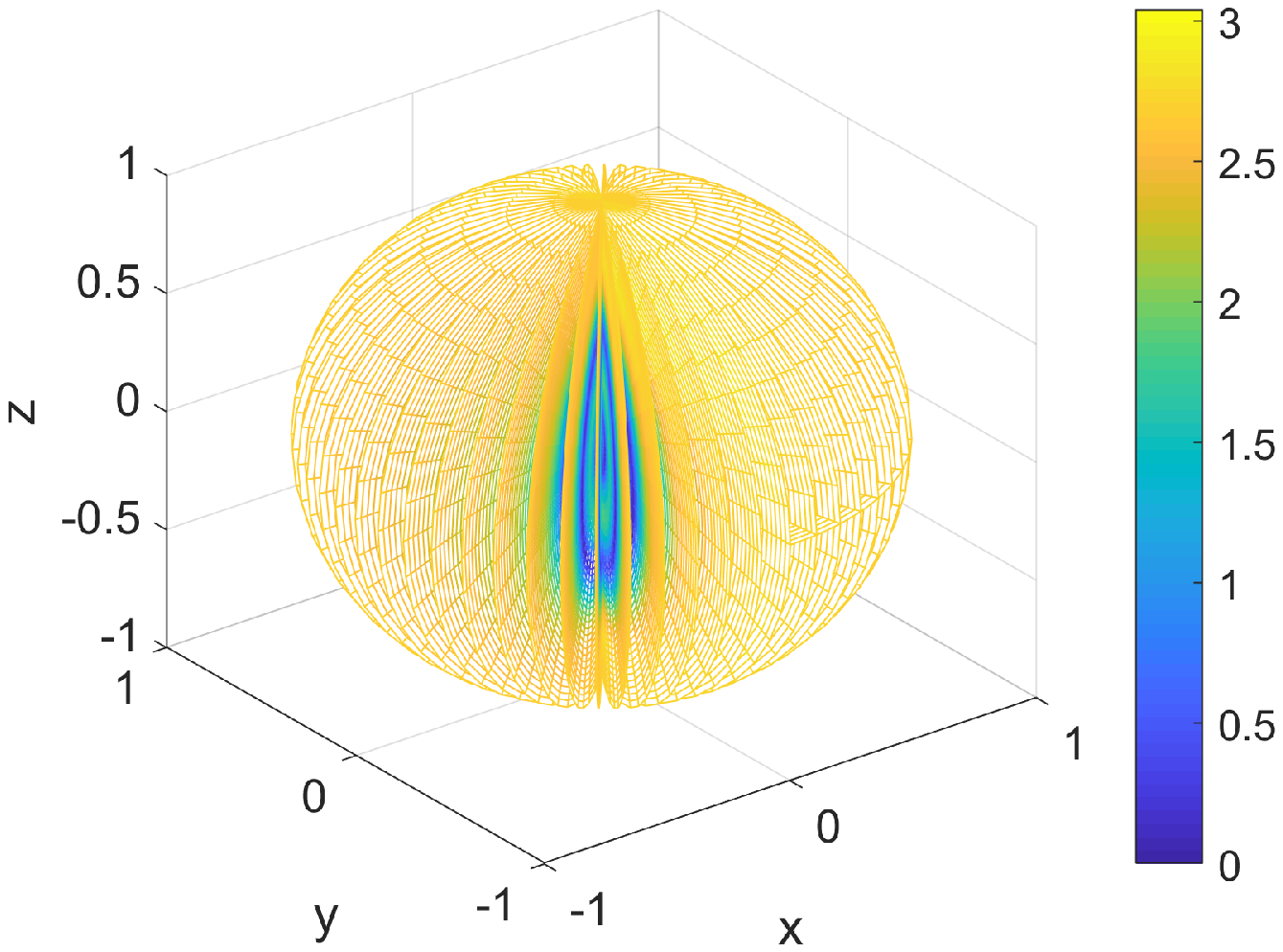}
		&\includegraphics[width=0.5\textwidth]{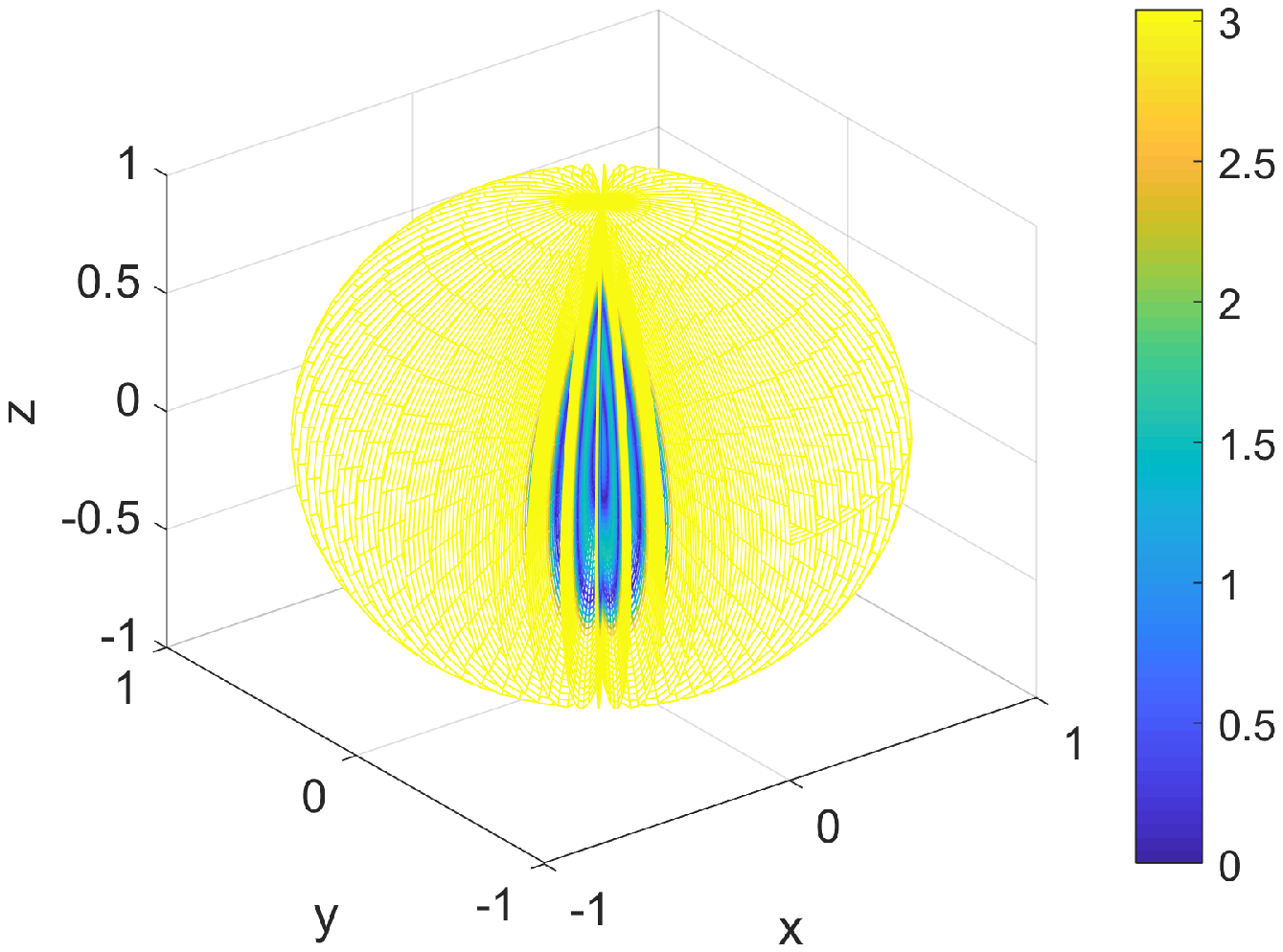}		\\
		(a) t=0.1s & (b) t=0.2s 
	\end{tabular}\\
	\begin{tabular}{@{}c@{}c@{}}
		\includegraphics[width=0.5\textwidth]{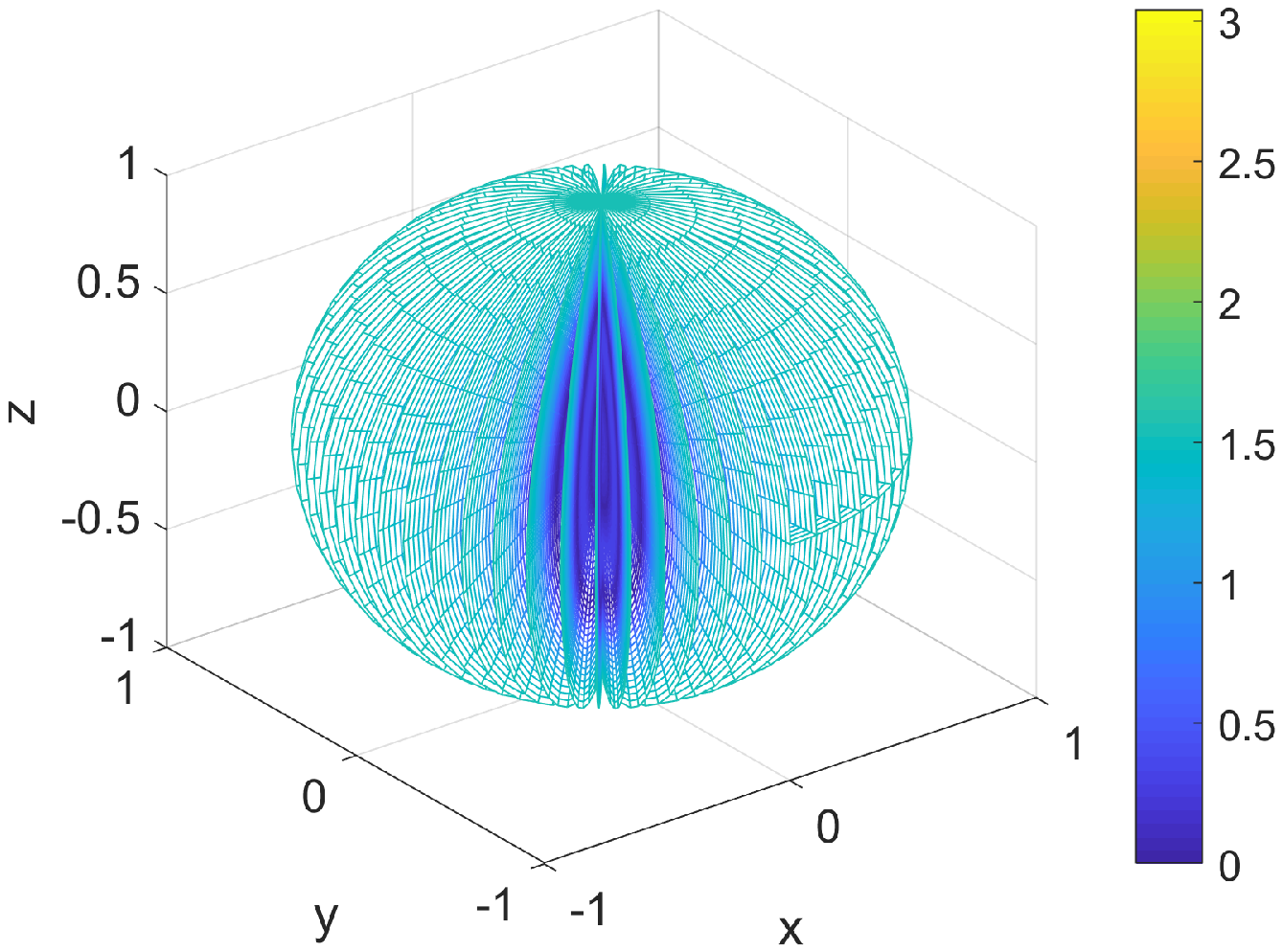}
		&\includegraphics[width=0.5\textwidth]{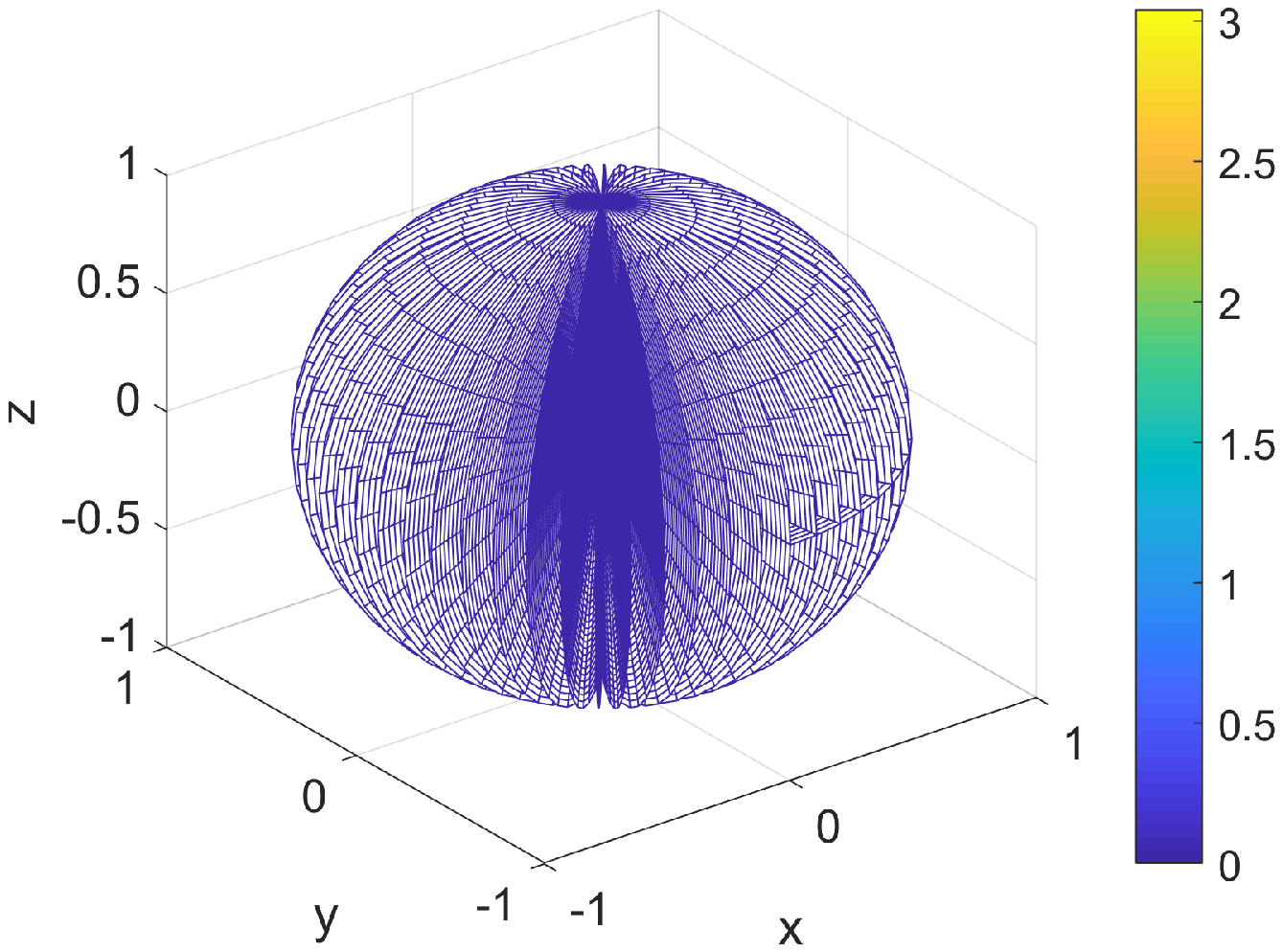}		\\ 
		(c) t=0.4s & (d) t=2s
	\end{tabular}\\
\begin{tabular}{@{}c@{}c@{}}
	\includegraphics[width=0.5\textwidth]{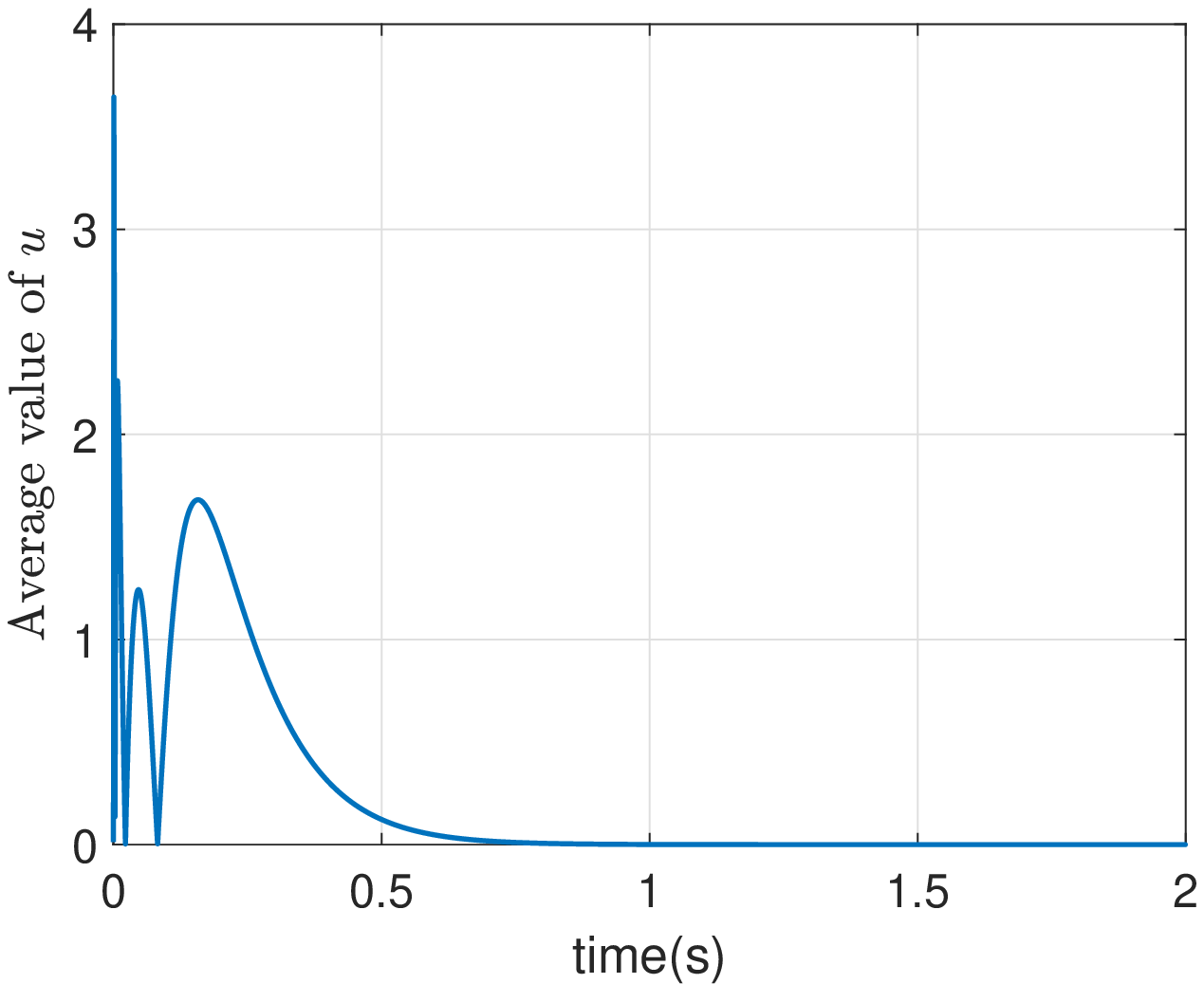}
	&\includegraphics[width=0.5\textwidth]{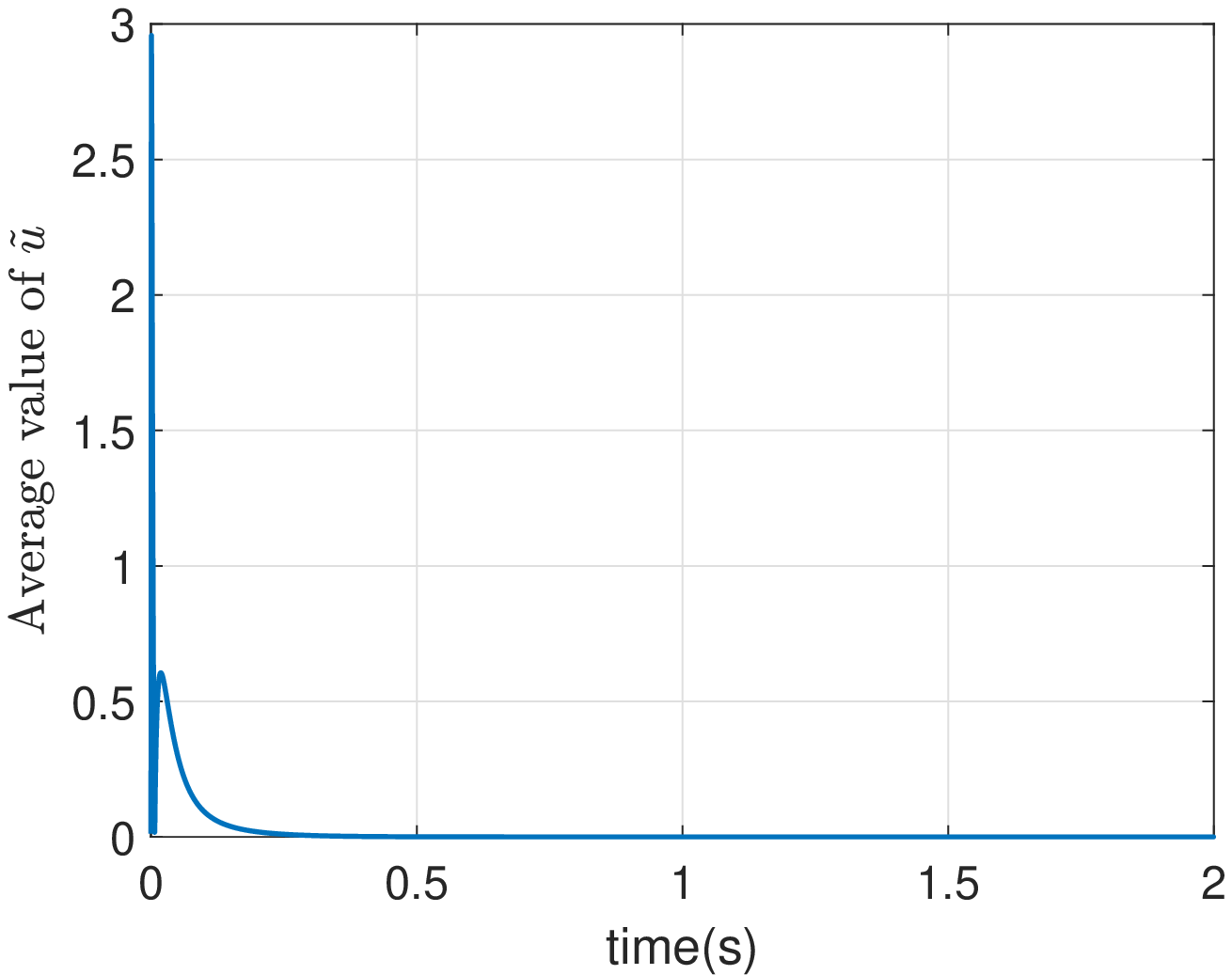}		\\
	(e) Mean (closed-loop system)   & (f)  Mean (observer error $\tilde u$)
\end{tabular}\\
	\caption{Closed-loop evolution using output feedback control.
		(a)-(d) Transient states. Note that the different
		upper ranges of the color bars in Fig. \ref{fig_openloop} and \ref{fig_close}. (e)
		Mean of system state u. (f) Mean of observer error $\tilde u$.}
	\label{fig_close}
\end{figure}

\begin{figure} 
	\centering
	\footnotesize
	\begin{tabular}{@{}c@{}c@{}}
		\includegraphics[width=0.5\textwidth]{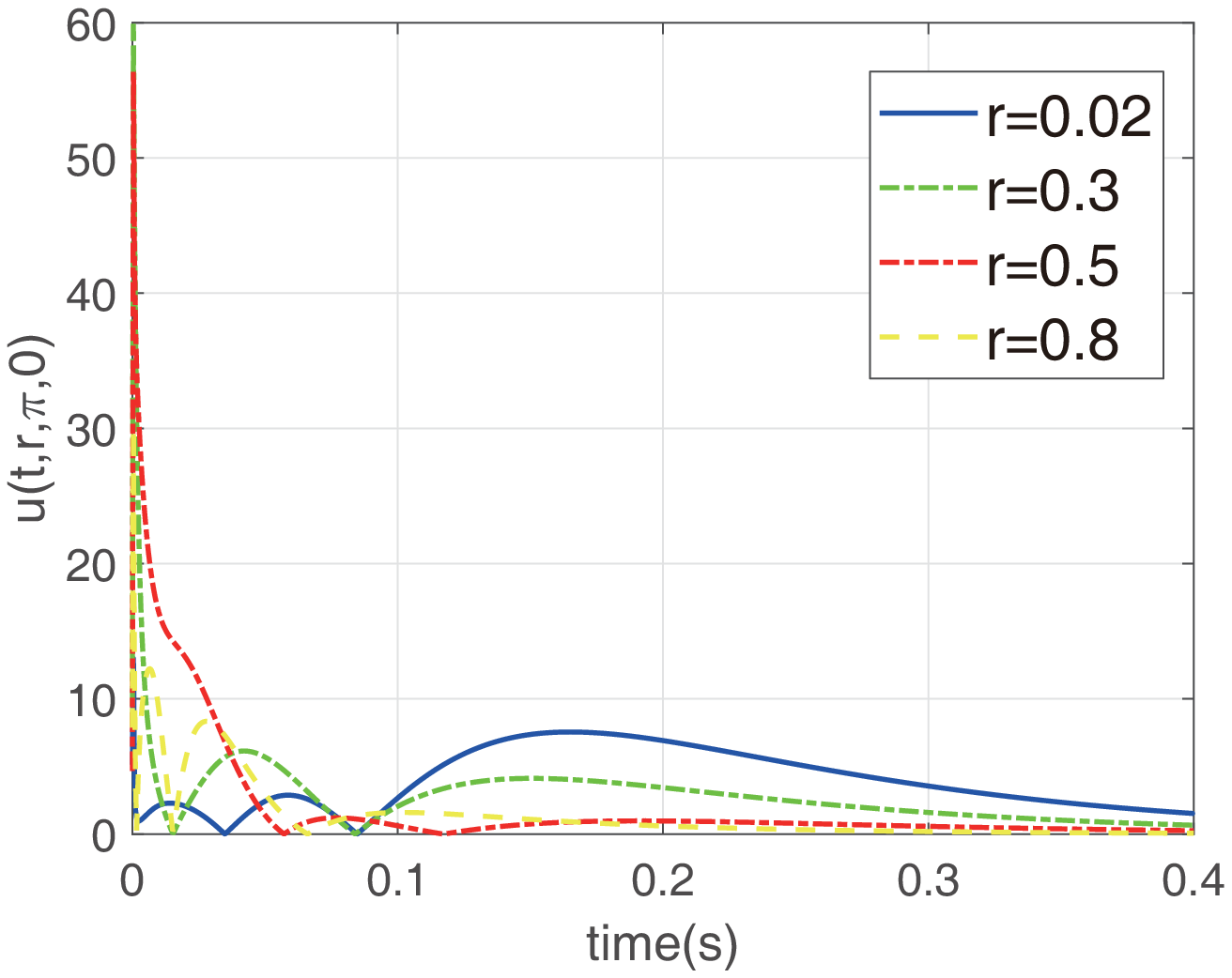}
		&\includegraphics[width=0.5\textwidth]{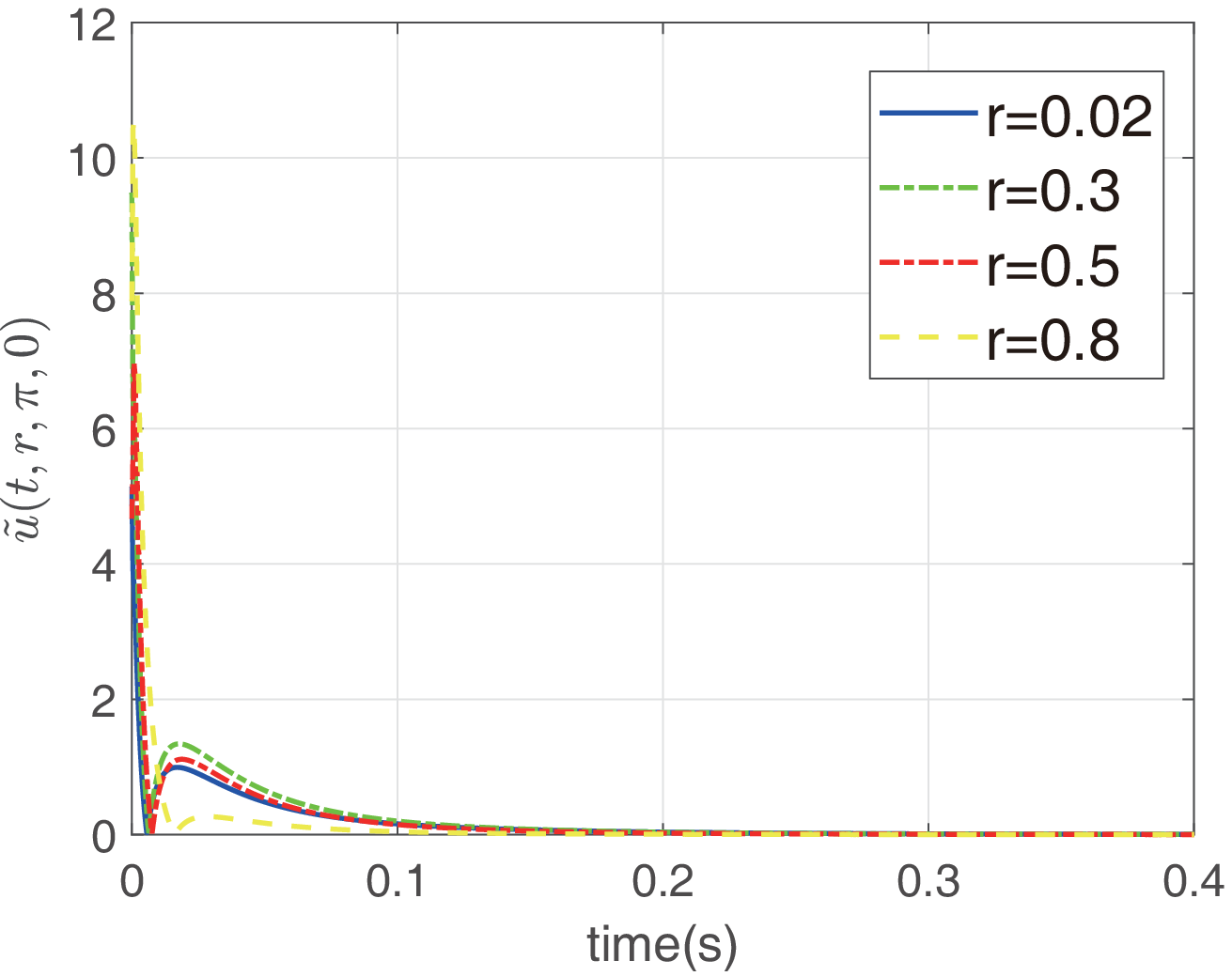}		\\
		(a) Actual states at $(r,\pi,0)$   & (b) Observer errors at $(r,\pi,0)$
	\end{tabular}\\
	\begin{tabular}{@{}c@{}c@{}}
		\includegraphics[width=0.5\textwidth]{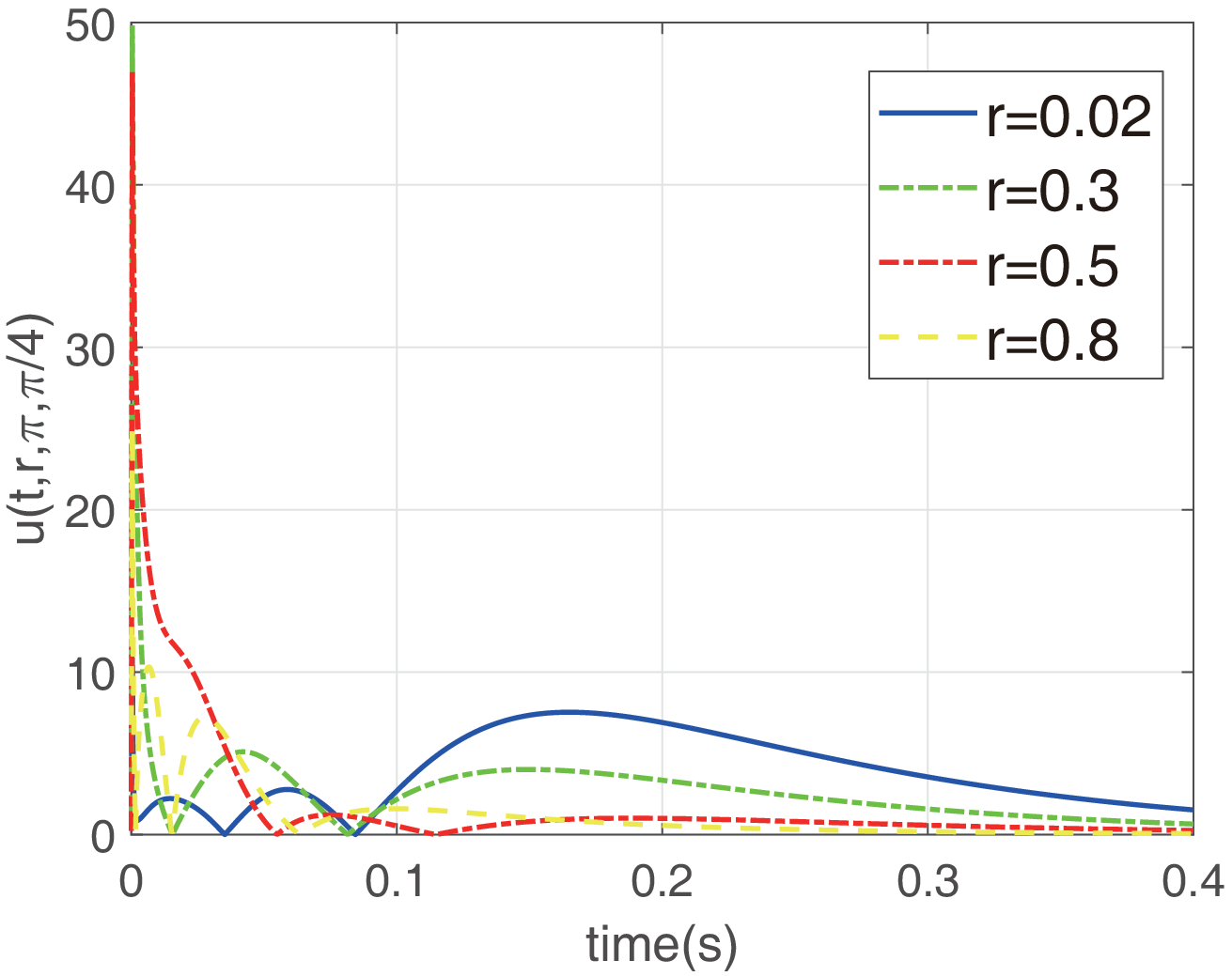}
		&\includegraphics[width=0.5\textwidth]{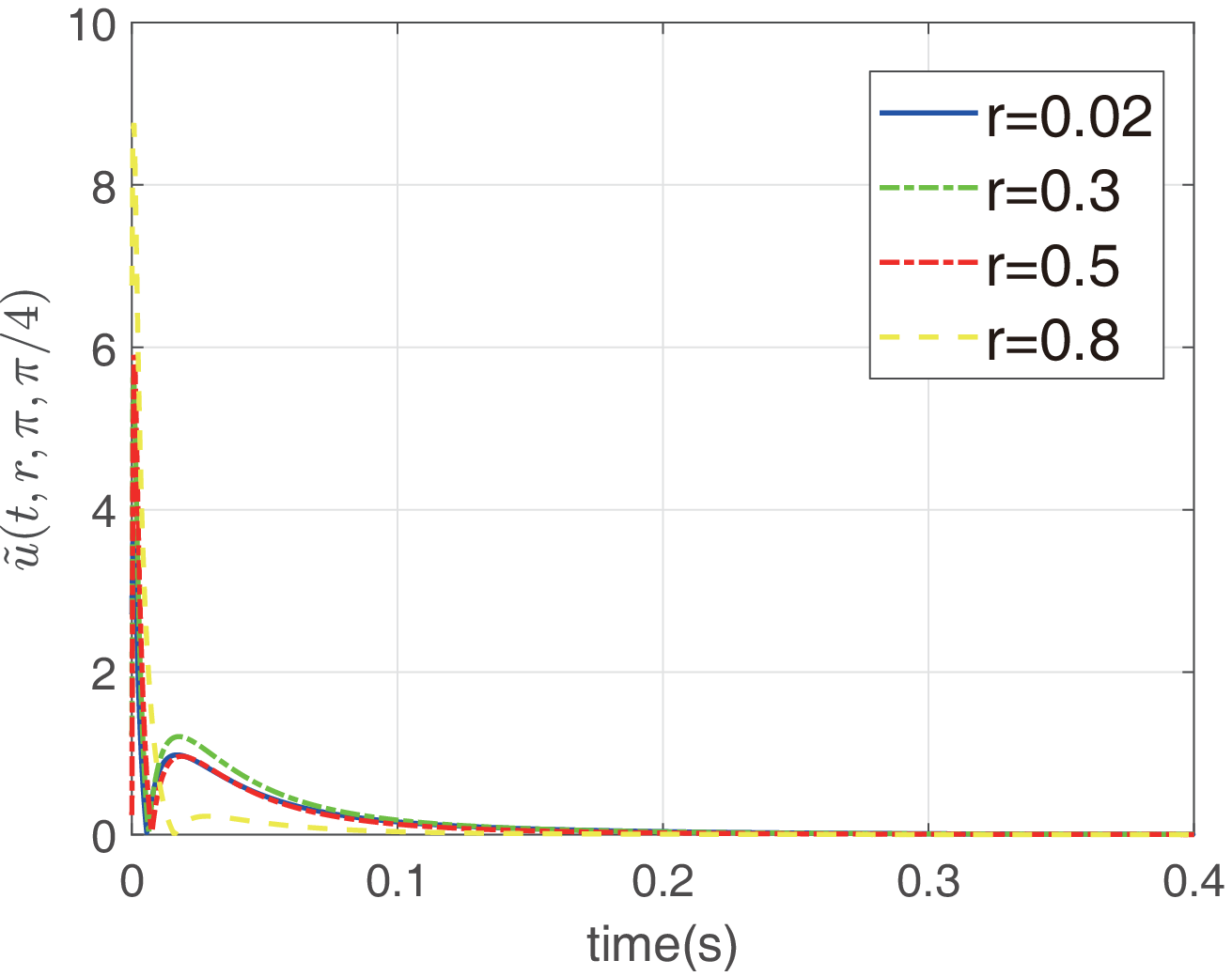}		\\
		(c) Actual states at $(r,\pi,\pi/4)$   & (d) Observer errors at $(r,\pi,\pi/4)$ 
	\end{tabular}\\
	\caption{The details of closed-loop evolution at different $r$ or $\vec \theta$. (a) (c) Actual states. (b) (d) Observer errors between the actual and estimated states.}
	\label{fig_detail}
\end{figure}

\begin{figure} 
	\centering
	\footnotesize
	\begin{tabular}{@{}c@{}c@{}}
		\includegraphics[width=0.5\textwidth]{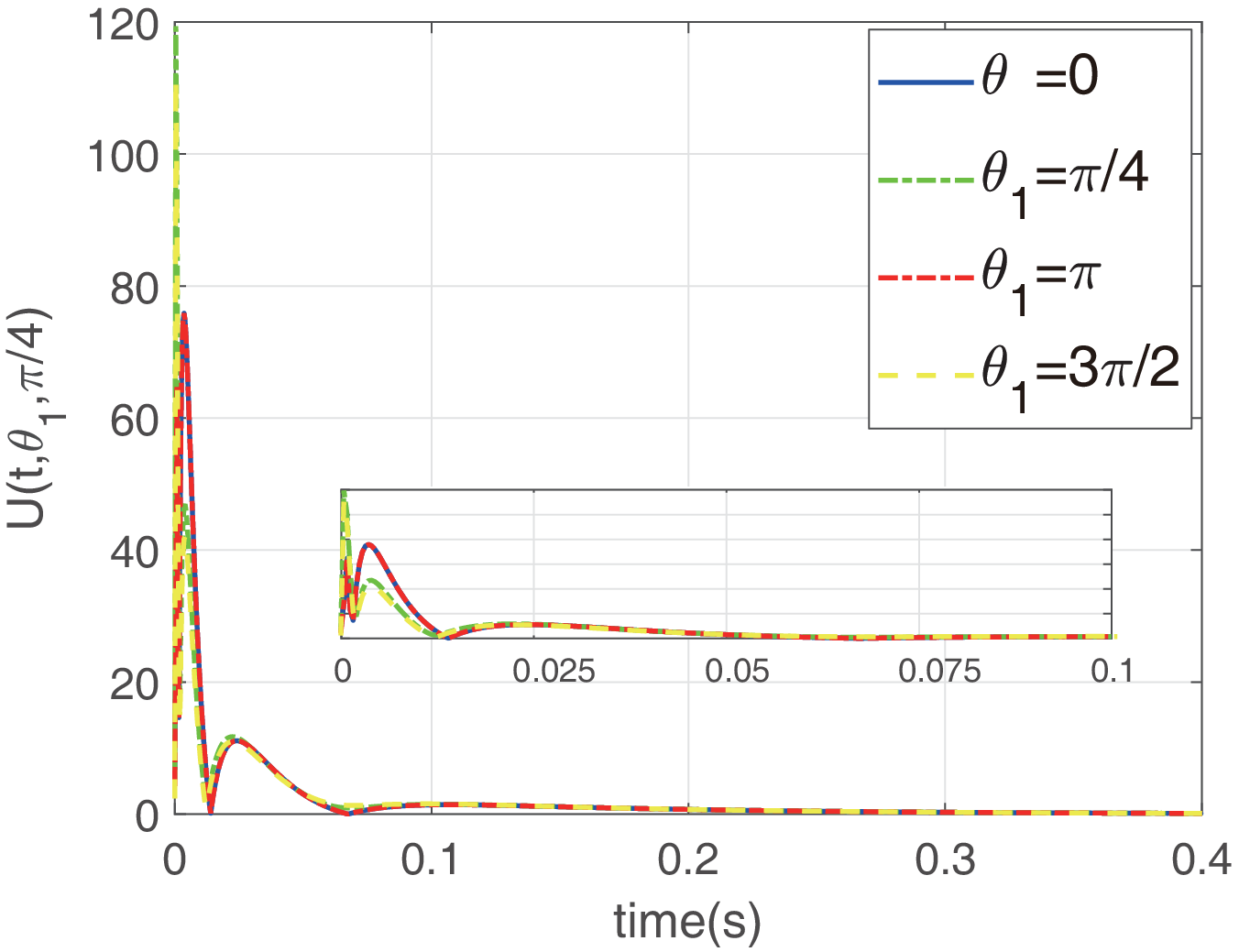}
		&\includegraphics[width=0.5\textwidth]{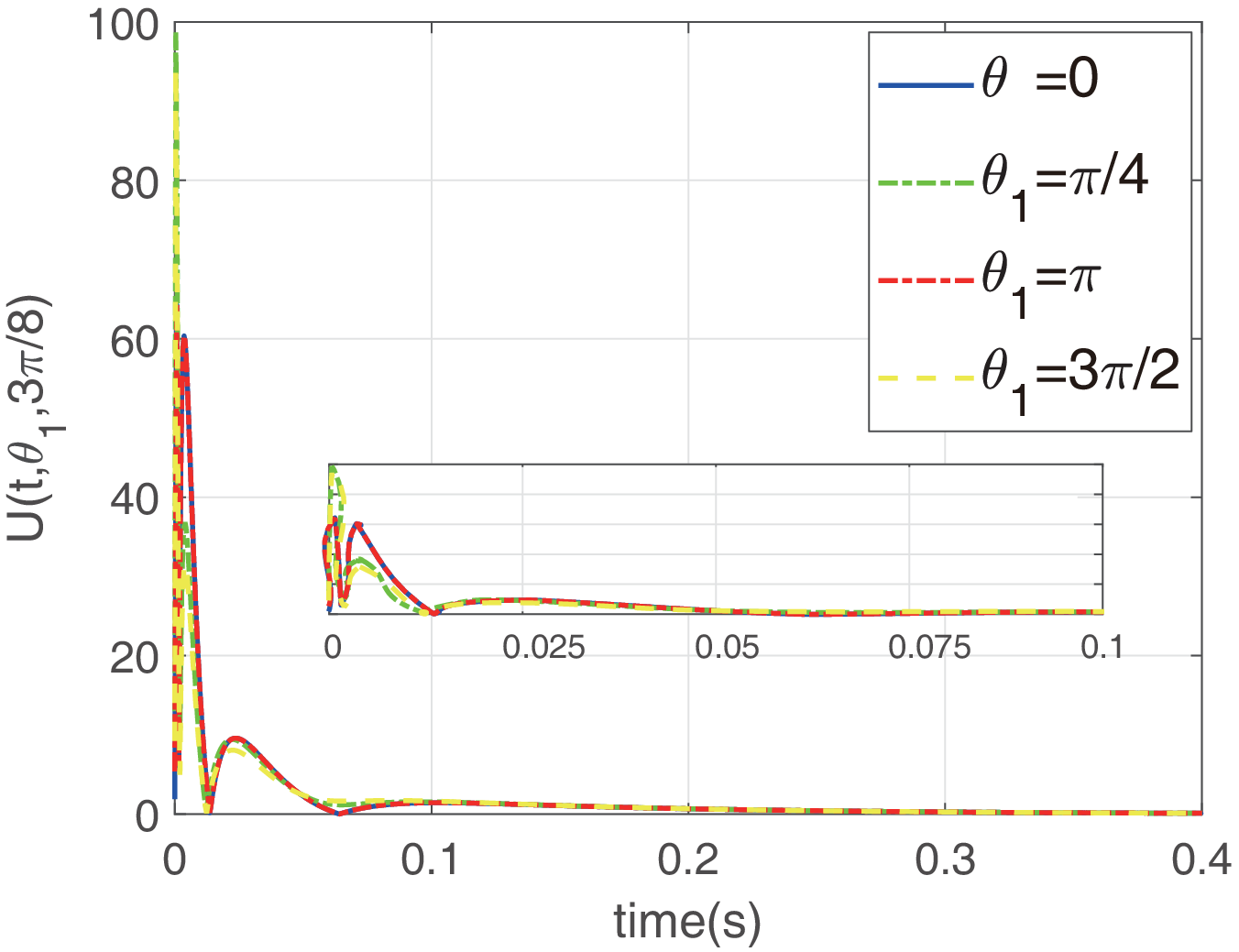}		\\
		(a) Control effort at $(\theta_1, \pi/4)$   & (b)  Control effort at $(\theta_1, 3\pi/8)$ 
	\end{tabular}\\
	\caption{The control effort at different $\vec \theta$.}
	\label{fig_control}
\end{figure}

Fig.~\ref{fig_openloop} and Fig.~\ref{fig_close} illustrate the transients of open-loop and closed-loop responses at different times, respectively, where the colour denotes the value of the position at this time. The evolutions of average value of $u$ are plotted in Fig.~\ref{fig_openloop}(d) and Fig.~\ref{fig_close}(e), respectively. (Note that the ranges of color bars are different. And thus avoid the appearance
of almost similar colors in fig. \ref{fig_close} of using same upper limit.)
Comparing the open-loop and closed-loop  evolution, the validity of proposed method is illustrated more intuitively. Fig.~\ref{fig_close}(f) shows the average of the observation errors, from which it can be found that the system begins to converge to its zero equilibrium after the observation error has already settled to zero as well.
The evolutions at different layers, namely $r = 0.002$, $r = 0.3$,  $r = 0.5$ and $r = 0.8$, are shown in Fig.~\ref{fig_detail} (a), (c) ,
as well as the observer errors  are presented in Fig.~\ref{fig_detail} (b), (d). For clarity, only the first $0.4$s of response are shown here.  Fig.~\ref{fig_control} depicts control effort at the boundary. It can be seen that the system driven by the proposed boundary control  eventually  converges after a short-term fluctuation.

\section{Conclusion}\label{sec-conclusions}
We have shown a  design to stabilize a radially-varying reaction-diffusion equation on an $n$-ball, by using an output-feedback boundary control law (with boundary measurements as well) designed through a backstepping method.  The radially-varying case proves to be a challenge as the kernel equations become singular in the radius; when applying the backstepping method, the same type of singularity appears in the kernel equations and successive approximations become difficult to use. Using a power series approach, a solution is found, thus providing a numerical method that can be readily applied, to both  control and observer boundary design. In addition the required conditions for the radially-varying coefficients are revealed (analyticity and evenness). 

In practice, this result can be of interest for deployment of multi-agent systems, by following the spirit of~\cite{jie}; thus, the radial domain mirrors a radial topology of interconnected agents which follow the reaction-diffusion dynamics to converge to equilibria, that represent different deployment profiles. Since one can choose the plant as desired (thus setting the behaviour of the agents), using analytic reaction coefficients is not actually a restriction, but rather opens the door to richer families of deployment profiles compared with the constant-coefficient case of~\cite{jie}.

On the other hand, the theoretical side of the result needs to be further investigated; an avenue of research that can be explored is the relaxation of the analyticity hypothesis by using reaction coefficients belonging to the Gevrey family; the kernels can then be analyzed to verify if they are still analytic, or rather Gevrey-type kernels, or simply do not converge.



\section*{References}
\begin{thebibliography}{xx}

\bibitem{abramowitz}
M.~Abramowitz and I.~A.~Stegun, {\em Handbook of mathematical functions},  9th Edition, Dover, 1965.

%
%
%
%


\bibitem{rijke}
G. Andrade, R. Vazquez and D. Pagano, ``Backstepping stabilization of a linearized ODE-PDE Rijke tube model,'' {\em Automatica}, vol. 96, 98-109, 2018.

\bibitem{ascencio}
P. Ascencio, A. Astolfi and T. Parisini, ``Backstepping PDE design: A convex optimization approach,'' {\em IEEE Transactions on Automatic Control}, vol. 63, pp. 1943--1958, 2018.

\bibitem{auriol}J. Auriol and F. Di Meglio, ``Minimum time control of heterodirectional linear coupled hyperbolic PDEs,'' {\em Automatica}, vol. 71, pp. 300--307, 2016.


\bibitem{atkinson}
K.Atkinson and W.~Han, {\em Spherical Harmonics and Approximations on the Unit Sphere: An Introduction}, Springer, 2012.

\bibitem{Barbu}
V.~Barbu, ``Boundary stabilization of equilibrium solutions to parabolic equations,'' {\em IEEE Transactions on Automatic Control}, vol. 58, pp. 2416--2420, 2013.

\bibitem{butkov}
E.~Butkov, {\em Mathematical physics}, Addison-Wesley, 1995.

\bibitem{leo}
L.~Camacho-Solorio, R.~Vazquez, M.~Krstic, ``Boundary observers for coupled diffusion-reaction systems with prescribed convergence rate,'' accepted in {\em Systems and Control Letters}, 2019.
%
%
%
%
%
%
%
%

\bibitem{howell}
K.~B.~Howell,  {\em Ordinary Differential Equations: An Introduction to the Fundamentals}, CRC Press, 2nd edition, 2019.

\bibitem{knuth}
D.~E.~Knuth, ``Two notes on notation,'' {\em The American Mathematical Monthly}, vol. 99(5), pp. 403--422, 1992.


\bibitem{krstic}
M.~Krstic and A.~Smyshlyaev, {\em Boundary Control of PDEs},  SIAM, 2008.

\bibitem{krstic5}
M.~Krstic, {\em Delay Compensation for Nonlinear, Adaptive, and PDE Systems}, Birkhauser, 2009.

\bibitem{long-nonlinear}
L.~Hu, R.~Vazquez, F.~Di~Meglio, and M.~Krstic, ``Boundary exponential stabilization of 1-D inhomogeneous quasilinear hyperbolic systems,'' {\em SIAM J. Control Optim.}, vol. 57(2), pp. 963--998, 2019.

\bibitem{meurer}
T.~Meurer and M.~Krstic, ``Finite-time multi-agent deployment: A nonlinear PDE motion planning approach,'' {\em Automatica}, vol. 47, pp. 2534--2542, 2011.

\bibitem{meurer2}
T.~Meurer, {\em Control of Higher-Dimensional PDEs: Flatness and Backstepping Designs}, Springer, 2013.

\bibitem{scott}
S.J.~Moura, n.A.~Chaturvedi, and M.~Krstic, ``PDE estimation techniques for advanced battery management systems---Part I: SOC estimation,''  {\em Proceedings of the 2012 American Control Conference}, 2012.


\bibitem{jie}
J.~Qi, R.~Vazquez and M.~Krstic, ``Multi-Agent deployment in 3-D via PDE control,'' {\em  IEEE Transactions on Automatic Control}, vol. 60 (4), pp. 891--906, 2015.

\bibitem{jie2}
J.~Qi, M.~Krstic and S. Wang, ``Stabilization of reaction-diffusion PDE distributed actuation and input delay,'' {\em Proceedings of the 2018 IEEE Conference on
Decision and Control (CDC)}, 2018.


\bibitem{triggiani}
R.~Triggiani, ``Boundary feedback stabilization of parabolic equations.{\em Appl. Math. Optimiz.},  vol. 6, pp. 201--220, 1980.

\bibitem{vazquez}
R.~Vazquez and M.~Krstic, {\em Control of Turbulent and Magnetohydrodynamic Channel Flow}.  Birkhauser, 2008.

\bibitem{vazquez2}
R.~Vazquez and M.~Krstic, ``Control of 1-D parabolic PDEs with Volterra nonlinearities --- Part I: Design,'' {\em  Automatica}, vol. 44, pp. 2778--2790, 2008.

\bibitem{Prieur}
F.~Bribiesca~Argomedo, C.~Prieur, E.~Witrant, and S.~Bremond, ``A strict control Lyapunov function for a diffusion equation with time-varying distributed coefficients,'' {\em IEEE Trans. Autom. Control}, vol. 58, pp. 290--303, 2013.

\bibitem{convloop}
R.~Vazquez and M.~Krstic, ``Boundary observer for output-feedback stabilization of thermal convection loop,'' {\em 
IEEE Trans. Control Syst. Technol.}, vol.18, pp. 789--797, 2010.

\bibitem{nball}
R. Vazquez and M. Krstic, ``Boundary control of reaction-diffusion PDEs on balls in spaces of arbitrary dimensions,'' {\em ESAIM:Control, Optimization and Calculus of Variations}, vol. 22, No. 4, pp. 1078--1096, 2016.

\bibitem{sphere}
R. Vazquez and M. Krstic, ``Boundary control and estimation of reaction-diffusion equations on the sphere under revolution symmetry conditions,'' {\em International Journal of Control}, vol. 92(1), pp. 2--11, 2019.

\bibitem{catalan}
R. Vazquez and M. Krstic, ``Boundary control of a singular reaction-diffusion equation on a disk,'' {\em CPDE 2016 (2nd IFAC Workshop on Control of Systems Governed by Partial Differential Equations)}, 2016.

 \bibitem{vazquez-coron}
R.~Vazquez, E.~Trelat and J.-M.~Coron, ``Control for fast and stable laminar-to-high-Reynolds-numbers transfer in a 2D navier-Stokes channel flow,'' {\em Disc. Cont. Dyn. Syst. Ser. B}, vol. 10, pp. 925--956, 2008.

 \bibitem{vazquez-zhang}
R.~Vazquez, M.~Krstic, J.~Zhang and J.~Qi, ``Stabilization of a 2-D reaction-diffusion equation with a coupled PDE evolving on its boundary,'' {\em Proceedings of the 2019 IEEE Conference on
Decision and Control (CDC)}, 2019.

\bibitem{jing-paper} R.~Vazquez, M.~Krstic, J.~Zhang and J.~Qi,  "Output Feedback Control of Radially-Dependent Reaction-Diffusion PDEs on Balls of Arbitrary Dimensions." IFAC-PapersOnLine 53, no. 2,pp. 7635--7640, 2020.

\end{thebibliography}

\end{document}